\documentclass[leqno,12pt]{amsart}
\usepackage{amssymb} 
\theoremstyle{plain} 
\newtheorem{theorem}{\indent\sc Theorem}[section] 
\newtheorem{lemma}[theorem]{\indent\sc Lemma}

\newtheorem{proposition}[theorem]{\indent\sc Proposition}

\theoremstyle{definition} 

\newtheorem{remark}[theorem]{\indent\sc Remark}

\newtheorem{problem}[theorem]{\indent\sc Problem}
\newtheorem{assumption}[theorem]{\indent\sc Assumption}
\begin{document}

\title[Coupled Painlev\'e VI systems]{Coupled Painlev\'e VI systems in dimension four\\
 with affine Weyl group symmetry\\
 of types $B_6^{(1)}$, $D_6^{(1)}$ and $D_7^{(2)}$ \\}
\author{Yusuke Sasano }

\renewcommand{\thefootnote}{\fnsymbol{footnote}}
\footnote[0]{2000\textit{ Mathematics Subjet Classification}.
 34M55, 34M45, 58F05, 32S65, 14E05, 20F55.}

\maketitle

\begin{abstract}
We find four kinds of six-parameter family of coupled Painlev\'e VI systems in dimension four with affine Weyl group symmetry of types  $B_6^{(1)}$, $D_6^{(1)}$ and $D_7^{(2)}$. Each system is the first example which gave higher-order Painlev\'e equations of types $B_l^{(1)},D_l^{(1)}$ and $D_l^{(2)}$, respectively.  Each system can be expressed as a polynomial Hamiltonian system.  We show that these systems are equivalent by an explicit birational and symplectic transformation, respectively. By giving each holomorphy condition, we can recover each system. These symmetries, holomorphy conditions and invariant divisors are new. We also give an explicit description of a confluence process from the system of type $D_6^{(1)}$ to the system of type $A_5^{(1)}$ by taking the coupling confluence process from the Painlev\'e VI system to the Painlev\'e V system.
\end{abstract}

\section{Introduction}
In 1912, considering the significant problem of searching for higher-order analogues of the Painlev\'e equations, Garnier discovered  a series of systems of nonlinear partial differential equations, which can be considered as a generalization of the Painlev\'e VI equation from the viewpoint of monodromy preserving deformations of the second-order linear ordinary differential equations, now called the Garnier system (see \cite{Garnier}).

From the viewpoint of affine Weyl groups, a series of systems of nonlinear ordinary differential equations with affine Weyl group symmetry of type $A_l^{(1)}$ were studied (cf. \cite{Adler}). This series gives a generalization of Painlev\'e equations $P_{IV}$ and $P_{V}$ to higher orders.

The Painlev\'e VI equation has symmetry under the affine Weyl group of type $D_4^{(1)}$. On the other hand, the generalizations of the systems of type $A_l^{(1)}$ do not include the Painlev\'e VI equation. Thus, it is an important remaining problem to find a generalization of the Painlev\'e VI equation for which the symmetries can be established. In the present paper, we find a 6-parameter family of coupled Painlev\'e VI systems in dimension four with affine Weyl group symmetry of type $D_6^{(1)}$. Our differential system is equivalent to a Hamiltonian system given by
\begin{equation}\label{1}
\frac{dx}{dt}=\frac{\partial H}{\partial y}, \ \ \frac{dy}{dt}=-\frac{\partial H}{\partial x}, \ \ \frac{dz}{dt}=\frac{\partial H}{\partial w}, \ \ \frac{dw}{dt}=-\frac{\partial H}{\partial z}
\end{equation}
with the polynomial Hamiltonian
\begin{align}\label{2}
\begin{split}
H &=H_{VI}(x,y,t;\alpha_0,\alpha_1,\alpha_2,\alpha_3,\alpha_4)+H_{VI}(z,w,t;\beta_0,\beta_1,\beta_2,\beta_3,\beta_4)\\
  &+\frac{2(x-t)yz\{(z-1)w+\beta_2\}}{t(t-1)}.
\end{split}
\end{align}
Here $x,y,z,w$ denote unknown complex variables, and $\alpha_0,\alpha_1,\alpha_2,\alpha_4,\beta_2,\beta_3,\beta_4$ are complex parameters satisfying the relation:
\begin{equation}\label{3}
\begin{aligned}
\alpha_0+\alpha_1+2\alpha_2+2(\alpha_4-\beta_4)+2\beta_2+\beta_3+\beta_4=1.
\end{aligned}
\end{equation}
We remark that for this system we tried to seek its first integrals of polynomial type with respect to $x,y,z,w$. However, we can not find. Of course, the Hamiltonian $H$ is not its first integral.

\begin{remark}
{\rm
The parameters $\alpha_0,\alpha_1,\dots ,\alpha_4,\beta_0,\beta_1,\dots ,\beta_4$ satisfy the relations{\rm : \rm}
\begin{align}\label{4}
\begin{split}
&\alpha_0+\alpha_1+2\alpha_2+\alpha_3+\alpha_4=1, \  \beta_0+\beta_1+2\beta_2+\beta_3+\beta_4=1,
\end{split}\\
\begin{split}\label{44}
&\alpha_1+2\alpha_2+\alpha_4-\beta_1-\beta_4=0, \ \alpha_3-\alpha_4-2\beta_2-\beta_3+\beta_4=0.
\end{split}
\end{align}
From these relations, it is easy to see that the parameters $\alpha_3,\alpha_4,\beta_0,\beta_1$ are described as linear combinations of the basis elements $\alpha_0,\alpha_1,\alpha_2,\beta_2,\beta_3,\beta_4$, which is explicitly given as follows{\rm : \rm}
\begin{align}\label{5}
\begin{split}
&\alpha_3=\frac{1-\alpha_0-\alpha_1-2\alpha_2+2\beta_2+\beta_3-\beta_4}{2}, \ \alpha_4=\frac{1-\alpha_0-\alpha_1-2\alpha_2-2\beta_2-\beta_3+\beta_4}{2},\\
&\beta_0=\frac{1+\alpha_0-\alpha_1-2\alpha_2-2\beta_2-\beta_3-\beta_4}{2}, \ \beta_1=\frac{1-\alpha_0+\alpha_1+2\alpha_2-2\beta_2-\beta_3-\beta_4}{2}.
\end{split}
\end{align}
\rm}
\end{remark}
The relations \eqref{4} are well-known as the parameter's relation of the Painlev\'e VI system, and the relations \eqref{44} are new. This representation of type $D_6^{(1)}$ can be constructed by coupling two copies of the $D_4^{(1)}$ root system by the relations \eqref{4} and \eqref{44} (see Theorem \ref{1.1}). This representation is new.

The symbol $H_{VI}(q,p,t;\delta_0,\delta_1,\delta_2,\delta_3,\delta_4)$ denotes the Hamiltonian of the second-order Painlev\'e VI equation given by
\begin{align}\label{6}
\begin{split}
&H_{VI}(q,p,t;\delta_0,\delta_1,\delta_2,\delta_3,\delta_4)\\
&=\frac{1}{t(t-1)}[p^2(q-t)(q-1)q-\{(\delta_0-1)(q-1)q+\delta_3(q-t)q\\
&+\delta_4(q-t)(q-1)\}p+\delta_2(\delta_1+\delta_2)q]  \quad (\delta_0+\delta_1+2\delta_2+\delta_3+\delta_4=1). 
\end{split}
\end{align}

This system is the first example which gave higher-order Painlev\'e equations of type $D_l^{(1)}$.

We also give an explicit description of a confluence process to the system with affine Weyl group symmetry of type $A_5^{(1)}$ (cf. \cite{Adler}). Additionally our results here, we obtain a new approach to the study of various higher-order Painlev\'e equations, presented in a series of papers for which this is the first. These papers are aimed at a complete study of the following problem:
\begin{problem}
For each affine root system $A$ with affine Weyl group $W(A)$, find a system of differential equations for which $W(A)$ acts as its B{\"a}cklund transformations.
\end{problem}

We remark that the B{\"a}cklund transformations of the system of type $D_6^{(1)}$ satisfy
\begin{equation}\label{unvre}
s_i(g)=g+\frac{\alpha_i}{f_i}\{f_i,g\}+\frac{1}{2!} \left(\frac{\alpha_i}{f_i} \right)^2 \{f_i,\{f_i,g\} \}+\cdots \quad (g \in {\Bbb C}(t)[x,y,z,w]),
\end{equation}
where poisson bracket $\{,\}$ satisfies the relations:
\begin{equation}
\{y,x\}=\{w,z\}=1, \quad the \ others \ are \ 0.
\end{equation}
Since these B{\"a}cklund transformations have Lie theoretic origin, similarity reduction of a Drinfeld-Sokolov hierarchy admits such a B{\"a}cklund symmetry.

This paper is organized as follows. In Section 2, we state our motivation and main results. In Section 3, we present two types of a 6-parameter family of coupled Painlev\'e VI systems in dimension four with affine Weyl group symmetry of type $B_6^{(1)}$. In Section 4, we find a 6-parameter family of coupled Painlev\'e VI systems in dimension four with affine Weyl group symmetry of type $D_7^{(2)}$. In Section 5, we find an autonomous version of the system of type $D_6^{(1)}$. In Section 6, we will give a brief review on the systems of types $A_5^{(1)}$ and $A_4^{(1)}$. In Section 7, we will explain our approach for obtaining the system \eqref{1}. In Section 8, we will prove Theorems \ref{1.3} and \ref{1.4}.

\section{Motivation and main results}

In the works \cite{Sasa10,Sasa11,Sasa12}, the author studied higher-order Painlev\'e equations from a viewpoint of polynomial Hamiltonian systems. In the case of the second-order Painlev\'e equations, let us summarize the following important properties of the Painlev\'e equations.

\vspace{0.1cm}
{\bf Notation.}

$\bullet$ $H \in {\Bbb C}(t)[x,y]$, \qquad $\bullet$ $deg(H)$: degree with respect to $x,y$.

\begin{center}
\begin{tabular}{|c|c|c|c|c|c|} \hline 
symmetry & $W(D_4^{(1)})$  & $W(A_3^{(1)})$ & $W(A_2^{(1)})$  & $W(C_2^{(1)})$ & $W(A_1^{(1)})$    \\ \hline 
Painlev\'e equations & $P_{VI}$ & $P_{V}$ & $P_{IV}$ & $P_{III}$ & $P_{II}$    \\ \hline 
degree of Hamiltonian $H$ & 5  & 4 & 3 & 4 & 3 \\ \hline
\end{tabular}
\end{center}

We are interested in polynomial Hamiltonian systems and symmetry under the affine Weyl group, and wish to search for higher-order Painlev\'e systems with these favorable properties. As examples of higher-order Painlev\'e systems, Adler (see \cite{Adler2,Galina}) studied ordinary differential systems with affine Weyl group symmetry of type $A_l^{(1)}$. When $l=2$ (resp. 3), this system of type $A_2^{(1)}$ (resp. $A_3^{(1)}$) is equivalent to the fourth (resp. fifth) Painlev\'e equation $P_{IV}$ (resp. $P_{V}$). They are considered to be higher-order versions of $P_V$ (resp. $P_{IV}$) when $l$ is odd (resp. even). These two examples motivated the author to find examples of higher-order versions other than $P_V$ and $P_{IV}$.  At first, we study four-dimensional case. Let us summarize important properties of the system of types $A_5^{(1)}$ and $A_4^{(1)}$.
\begin{center}
\begin{tabular}{|c|c|c|} \hline 
symmetry & $W(A_5^{(1)})$  & $W(A_4^{(1)})$    \\ \hline 
Hamiltonian $H$ & $H_{V}(x,y,t)+H_{V}(z,w,t)$ & $H_{IV}(x,y,t)+H_{IV}(z,w,t)$    \\
         & $-2yzw+\frac{2xyzw}{t}$ & $+2yzw$    \\ \hline 
differential system & coupled Painlev\'e $V$ & coupled Painlev\'e $IV$     \\ \hline 
degree of Hamiltonian $H$ & 4  & 3 \\ \hline
\end{tabular}
\end{center}

These properties suggest the possibility that there exists a procedure for searching for such higher-order versions with symmetry under the affine Weyl group of type $D_6^{(1)}$. Here, let us consider the following problem.

\begin{problem}
Can we show existence of a coupled Painlev\'e VI system in dimension four satisfying the following assumptions $(A1),(A2)$? If yes, can we find it explicitly and is it unique?
\end{problem}

\begin{assumption}
$(A1)$ $deg(H)=5$ with respect to $x,y,z,w$.\\
$(A2)$ The system has symmetry under the affine Weyl group of type $D_6^{(1)}$.
\end{assumption}

To answer this, in this paper, we present a 6-parameter family of coupled Painlev\'e VI systems in dimension four with extended affine Weyl group symmetry of type $D_6^{(1)}$ explicitly given by
\begin{equation*}
  \left\{
  \begin{aligned}
   \frac{dx}{dt} &=\frac{\partial H}{\partial y}=\frac{1}{t(t-1)}\{2y(x-t)(x-1)x-(\alpha_0-1)(x-1)x-\alpha_3(x-t)x\\
                 & \qquad -\alpha_4(x-t)(x-1)+2(x-t)z((z-1)w+\beta_2)\},\\
   \frac{dy}{dt} &=-\frac{\partial H}{\partial x}=\frac{1}{t(t-1)}[-\{(x-t)(x-1)+(x-t)x+(x-1)x\}y^2+\{(\alpha_0-1)(2x-1)\\
                 & \qquad +\alpha_3(2x-t)+\alpha_4(2x-t-1)\}y-\alpha_2(\alpha_1+\alpha_2)-2yz((z-1)w+\beta_2)],\\
   \frac{dz}{dt} &=\frac{\partial H}{\partial w}=\frac{1}{t(t-1)}\{2w(z-t)(z-1)z-(\beta_0-1)(z-1)z-\beta_3(z-t)z\\
                 & \qquad -\beta_4(z-t)(z-1)+2(x-t)yz(z-1)\},\\
   \frac{dw}{dt} &=-\frac{\partial H}{\partial z}=\frac{1}{t(t-1)}[-\{(z-t)(z-1)+(z-t)z+(z-1)z\}w^2+\{(\beta_0-1)(2z-1)\\
                 &+\beta_3(2z-t)+\beta_4(2z-t-1)\}w-\beta_2(\beta_1+\beta_2)-2(x-t)y((2z-1)w+\beta_2)]
   \end{aligned}
  \right. 
\end{equation*}
with the polynomial Hamiltonian $H$ \eqref{2}.

\begin{theorem}\label{1.1}
The system \eqref{1} admits extended affine Weyl group symmetry of type $D_6^{(1)}$ as the group of its B{\"a}cklund transformations, whose generators $s_i \ (i=0,\dots ,6),{\pi}_j \ (j=1,\dots ,4)$ are explicitly given as follows{\rm : \rm}with the notations
$$
\gamma_1:=\alpha_4-\beta_4, \qquad  (*):=(x,y,z,w,t;\alpha_0,\alpha_1,\alpha_2,\gamma_1,\beta_2,\beta_3,\beta_4),
$$
\begin{align}
\begin{split}
        s_0: (*) &\rightarrow \left(x,y-\frac{\alpha_0}{x-t},z,w,t;-\alpha_0,\alpha_1,\alpha_2+\alpha_0,\gamma_1,\beta_2,\beta_3,\beta_4 \right),\\
        s_1: (*) &\rightarrow (x,y,z,w,t;\alpha_0,-\alpha_1,\alpha_2+\alpha_1,\gamma_1,\beta_2,\beta_3,\beta_4),\\
        s_2: (*) &\rightarrow \left(x+\frac{\alpha_2}{y},y,z,w,t;\alpha_0+\alpha_2,\alpha_1+\alpha_2,-\alpha_2,\gamma_1+\alpha_2,\beta_2,\beta_3,\beta_4 \right), \\
        s_3: (*) &\rightarrow \left(x,y-\frac{\gamma_1}{x-z},z,w+\frac{\gamma_1}{x-z},t;\alpha_0,\alpha_1,\alpha_2+\gamma_1,-\gamma_1,\beta_2+\gamma_1,\beta_3,\beta_4 \right), \\
        s_4: (*) &\rightarrow \left(x,y,z+\frac{\beta_2}{w},w,t;\alpha_0,\alpha_1,\alpha_2,\gamma_1+\beta_2,-\beta_2,\beta_3+\beta_2,\beta_4+\beta_2 \right), \\
        s_5: (*) &\rightarrow \left(x,y,z,w-\frac{\beta_3}{z-1},t;\alpha_0,\alpha_1,\alpha_2,\gamma_1,\beta_2+\beta_3,-\beta_3,\beta_4 \right), \\
        s_6: (*) &\rightarrow \left(x,y,z,w-\frac{\beta_4}{z},t;\alpha_0,\alpha_1,\alpha_2,\gamma_1,\beta_2+\beta_4,\beta_3,-\beta_4 \right),\\
        {\pi}_1: (*) &\rightarrow  (\frac{t(t-1)+t(x-t)}{x-t},-\frac{(x-t)((x-t)y+\alpha_2)}{t(t-1)},\frac{t(t-1)+t(z-t)}{z-t},\\
        &-\frac{(z-t)((z-t)w+\beta_2)}{t(t-1)},t;\alpha_1,\alpha_0,\alpha_2,\gamma_1,\beta_2,\beta_4,\beta_3),\\
        {\pi}_2: (*) &\rightarrow  \left(\frac{t}{z},-\frac{z(zw+\beta_2)}{t},\frac{t}{x},-\frac{x(xy+\alpha_2)}{t},t;\beta_3,\beta_4,\beta_2,\gamma_1,\alpha_2,\alpha_0,\alpha_1 \right), \\
        {\pi}_3: (*) &\rightarrow  (1-x,-y,1-z,-w,1-t;\alpha_0,\alpha_1,\alpha_2,\gamma_1,\beta_2,\beta_4,\beta_3), \\
        {\pi}_4: (*) &\rightarrow  (\frac{(t-1)x}{t-x},\frac{(t-x)(ty-xy-\alpha_2)}{t(t-1)},\frac{(t-1)z}{t-z},\frac{(t-z)(tw-zw-\beta_2)}{t(t-1)},\\
        &1-t;\alpha_1,\alpha_0,\alpha_2,\gamma_1,\beta_2,\beta_3,\beta_4).
        \end{split}
        \end{align}
\end{theorem}
We note that these transformations $s_i,\pi_j$ are birational and symplectic.

The B{\"a}cklund transformations of the system of type $D_6^{(1)}$ satisfy
\begin{equation}
s_i(g)=g+\frac{\alpha_i}{f_i}\{f_i,g\}+\frac{1}{2!} \left(\frac{\alpha_i}{f_i} \right)^2 \{f_i,\{f_i,g\} \}+\cdots \quad (g \in {\Bbb C}(t)[x,y,z,w]),
\end{equation}
where poisson bracket $\{,\}$ satisfies the relations:
\begin{equation}
\{y,x\}=\{w,z\}=1, \quad the \ others \ are \ 0.
\end{equation}
Since these B{\"a}cklund transformations have Lie theoretic origin, similarity reduction of a Drinfeld-Sokolov hierarchy admits such a B{\"a}cklund symmetry.

\begin{proposition}
The system \eqref{1} has the following invariant divisors\rm{:\rm}
\begin{center}
\begin{tabular}{|c|c|c|} \hline
parameter's relation & $f_i$ \\ \hline
$\alpha_0=0$ & $f_0:=x-t$  \\ \hline
$\alpha_1=0$ & $f_1:=x-\infty$  \\ \hline
$\alpha_2=0$ & $f_2:=y$  \\ \hline
$\alpha_4=\beta_4$ & $f_3:=x-z$  \\ \hline
$\beta_2=0$ & $f_4:=w$  \\ \hline
$\beta_3=0$ & $f_5:=z-1$  \\ \hline
$\beta_4=0$ & $f_6:=z$  \\ \hline
\end{tabular}
\end{center}
\end{proposition}
We note that when $\alpha_2=0$, we see that the system \eqref{1} admits a particular solution $y=0$, and when $\alpha_4=\beta_4$, after we make the birational and symplectic transformation:
\begin{equation}
x_3=x-z, \ y_3=y, \ z_3=z, \ w_3=w+y
\end{equation}
we see that the system \eqref{1} admits a particular solution $x_3=0$.

\begin{remark}
It is easy to see that the generators ${\pi}_2,{\pi}_3,{\pi}_4$ satisfy the relation{\rm : \rm}
\begin{equation*}
{\pi}_4={\pi}_2{\pi}_3{\pi}_2.
\end{equation*}
\end{remark}

\begin{remark}
Taking the coordinate system
\begin{equation*}
(X,Y,Z,W)=(1/x,-x(xy+\alpha_2),z,w),
\end{equation*}
it is easy to see that the transformation $s_1$ can be explicitly written as follows{\rm : \rm}
\begin{align*}
\begin{split}
&s_1: (X,Y,Z,W,t;\alpha_0,\alpha_1,\alpha_2,\gamma_1,\beta_2,\beta_3,\beta_4) \\
& \qquad \rightarrow (X,Y-\alpha_1/X,Z,W,t;\alpha_0,-\alpha_1,\alpha_1+\alpha_2,\gamma_1,\beta_2,\beta_3,\beta_4).
\end{split}
\end{align*}
\end{remark}

\begin{proposition}\label{P1.1}
Let us define the following translation operators
\begin{gather*}
\begin{gathered}
T_1:={\pi}_1s_5s_4s_3s_2s_1s_0s_2s_3s_4s_5, \quad T_2:=s_4s_6T_1s_6s_4, \quad T_3:=s_6T_1s_6,\\
T_4:={\pi}_2T_1{\pi}_2, \quad T_5:={\pi}_2T_2{\pi}_2, \quad T_6:={\pi}_2T_3{\pi}_2.
\end{gathered}
\end{gather*}
These translation operators act on parameters $\alpha_i,\gamma_1,\beta_j$ as follows{\rm : \rm}
\begin{align*}
T_1(\alpha_0,\alpha_1,\alpha_2,\gamma_1,\beta_2,\beta_3,\beta_4)&=(\alpha_0,\alpha_1,\alpha_2,\gamma_1,\beta_2,\beta_3,\beta_4)+(0,0,0,0,0,-1,1),\\
T_2(\alpha_0,\alpha_1,\alpha_2,\gamma_1,\beta_2,\beta_3,\beta_4)&=(\alpha_0,\alpha_1,\alpha_2,\gamma_1,\beta_2,\beta_3,\beta_4)+(0,0,0,1,-1,0,0),\\
T_3(\alpha_0,\alpha_1,\alpha_2,\gamma_1,\beta_2,\beta_3,\beta_4)&=(\alpha_0,\alpha_1,\alpha_2,\gamma_1,\beta_2,\beta_3,\beta_4)+(0,0,0,0,1,-1,-1),\\
T_4(\alpha_0,\alpha_1,\alpha_2,\gamma_1,\beta_2,\beta_3,\beta_4)&=(\alpha_0,\alpha_1,\alpha_2,\gamma_1,\beta_2,\beta_3,\beta_4)+(-1,1,0,0,0,0,0),\\
T_5(\alpha_0,\alpha_1,\alpha_2,\gamma_1,\beta_2,\beta_3,\beta_4)&=(\alpha_0,\alpha_1,\alpha_2,\gamma_1,\beta_2,\beta_3,\beta_4)+(0,0,-1,1,0,0,0),\\
T_6(\alpha_0,\alpha_1,\alpha_2,\gamma_1,\beta_2,\beta_3,\beta_4)&=(\alpha_0,\alpha_1,\alpha_2,\gamma_1,\beta_2,\beta_3,\beta_4)+(-1,-1,1,0,0,0,0).
\end{align*}
\end{proposition}

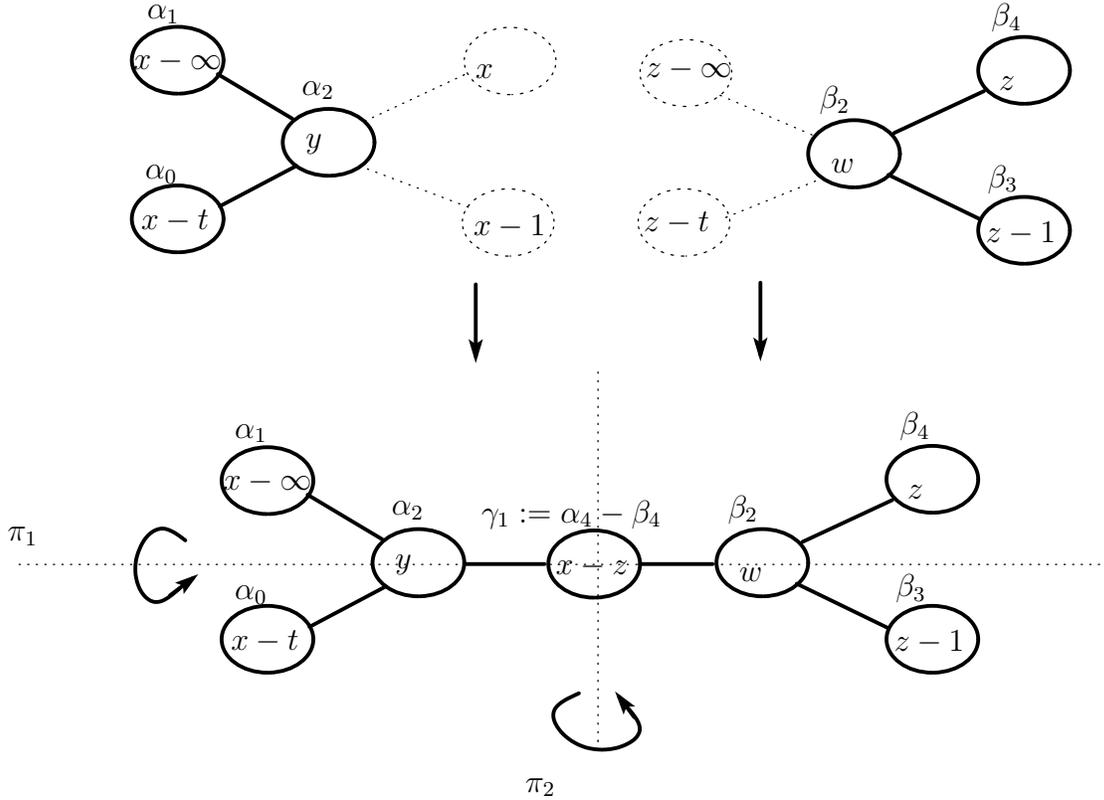
\begin{figure}[h]
\unitlength 0.1in
\begin{picture}(57.60,40.48)(12.60,-44.98)
%
\special{pn 20}%
\special{ar 2620 3014 240 175  1.5707963 6.2831853}%
\special{ar 2620 3014 240 175  0.0000000 1.5291537}%
%
\special{pn 20}%
\special{ar 2620 3845 240 175  1.5707963 6.2831853}%
\special{ar 2620 3845 240 175  0.0000000 1.5291537}%
%
\special{pn 20}%
\special{ar 3410 3444 240 176  1.5707963 6.2831853}%
\special{ar 3410 3444 240 176  0.0000000 1.5291537}%
%
\special{pn 20}%
\special{ar 5210 3444 240 176  1.5707963 6.2831853}%
\special{ar 5210 3444 240 176  0.0000000 1.5291537}%
%
\special{pn 20}%
\special{ar 6100 3007 240 175  1.5707963 6.2831853}%
\special{ar 6100 3007 240 175  0.0000000 1.5291537}%
%
\special{pn 20}%
\special{ar 6100 3845 240 175  1.5707963 6.2831853}%
\special{ar 6100 3845 240 175  0.0000000 1.5291537}%
%
\special{pn 20}%
\special{pa 2830 3088}%
\special{pa 3220 3321}%
\special{fp}%
%
\special{pn 20}%
\special{pa 2840 3780}%
\special{pa 3240 3568}%
\special{fp}%
%
\special{pn 20}%
\special{ar 4330 3451 240 176  1.5707963 6.2831853}%
\special{ar 4330 3451 240 176  0.0000000 1.5291537}%
%
\special{pn 20}%
\special{pa 3650 3451}%
\special{pa 4070 3451}%
\special{fp}%
%
\special{pn 20}%
\special{pa 4570 3451}%
\special{pa 4950 3451}%
\special{fp}%
%
\special{pn 20}%
\special{pa 5420 3335}%
\special{pa 5890 3116}%
\special{fp}%
%
\special{pn 20}%
\special{pa 5390 3554}%
\special{pa 5870 3787}%
\special{fp}%
%
\put(20.1000,-30.6500){\makebox(0,0)[lb]{}}%
\put(24.3000,-39.1200){\makebox(0,0)[lb]{$x-t$}}%
\put(32.9000,-35.0700){\makebox(0,0)[lb]{$y$}}%
\put(41.3000,-35.1700){\makebox(0,0)[lb]{$x-z$}}%
\put(50.9000,-35.4300){\makebox(0,0)[lb]{$w$}}%
\put(59.0000,-39.2100){\makebox(0,0)[lb]{$z-1$}}%
\put(59.7000,-31.1100){\makebox(0,0)[lb]{$z$}}%
\put(23.9000,-30.7500){\makebox(0,0)[lb]{$x-\infty$}}%
\put(39.7000,-46.6800){\makebox(0,0)[lb]{${\pi}_2$}}%
\put(12.6000,-33.5400){\makebox(0,0)[lb]{${\pi}_1$}}%
%
\special{pn 8}%
\special{pa 1320 3453}%
\special{pa 7020 3453}%
\special{dt 0.045}%
\special{pa 7020 3453}%
\special{pa 7019 3453}%
\special{dt 0.045}%
%
\special{pn 8}%
\special{pa 4350 2445}%
\special{pa 4350 4380}%
\special{dt 0.045}%
\special{pa 4350 4380}%
\special{pa 4350 4379}%
\special{dt 0.045}%
%
\special{pn 20}%
\special{pa 2190 3327}%
\special{pa 2161 3305}%
\special{pa 2133 3286}%
\special{pa 2104 3272}%
\special{pa 2077 3264}%
\special{pa 2049 3267}%
\special{pa 2023 3278}%
\special{pa 1998 3297}%
\special{pa 1976 3323}%
\special{pa 1958 3354}%
\special{pa 1943 3390}%
\special{pa 1932 3427}%
\special{pa 1928 3467}%
\special{pa 1929 3506}%
\special{pa 1937 3545}%
\special{pa 1950 3580}%
\special{pa 1969 3611}%
\special{pa 1992 3634}%
\special{pa 2017 3648}%
\special{pa 2044 3651}%
\special{pa 2071 3644}%
\special{pa 2099 3628}%
\special{pa 2127 3607}%
\special{pa 2156 3582}%
\special{pa 2160 3579}%
\special{sp}%
%
\special{pn 20}%
\special{pa 4250 4128}%
\special{pa 4216 4144}%
\special{pa 4183 4160}%
\special{pa 4155 4178}%
\special{pa 4133 4197}%
\special{pa 4119 4219}%
\special{pa 4115 4245}%
\special{pa 4120 4272}%
\special{pa 4134 4300}%
\special{pa 4154 4328}%
\special{pa 4180 4354}%
\special{pa 4212 4378}%
\special{pa 4248 4398}%
\special{pa 4288 4413}%
\special{pa 4329 4421}%
\special{pa 4372 4424}%
\special{pa 4414 4421}%
\special{pa 4454 4413}%
\special{pa 4490 4400}%
\special{pa 4521 4384}%
\special{pa 4546 4365}%
\special{pa 4563 4342}%
\special{pa 4570 4318}%
\special{pa 4566 4292}%
\special{pa 4553 4265}%
\special{pa 4533 4237}%
\special{pa 4510 4209}%
\special{pa 4510 4209}%
\special{sp}%
%
\special{pn 20}%
\special{pa 2110 3615}%
\special{pa 2230 3525}%
\special{fp}%
\special{sh 1}%
\special{pa 2230 3525}%
\special{pa 2165 3549}%
\special{pa 2187 3557}%
\special{pa 2189 3581}%
\special{pa 2230 3525}%
\special{fp}%
%
\special{pn 20}%
\special{pa 4540 4254}%
\special{pa 4460 4146}%
\special{fp}%
\special{sh 1}%
\special{pa 4460 4146}%
\special{pa 4484 4211}%
\special{pa 4492 4189}%
\special{pa 4516 4188}%
\special{pa 4460 4146}%
\special{fp}%
\put(24.5000,-36.6000){\makebox(0,0)[lb]{$\alpha_0$}}%
\put(24.5000,-28.1400){\makebox(0,0)[lb]{$\alpha_1$}}%
\put(32.7000,-32.1900){\makebox(0,0)[lb]{$\alpha_2$}}%
\put(50.3000,-32.3700){\makebox(0,0)[lb]{$\beta_2$}}%
\put(59.1000,-36.5100){\makebox(0,0)[lb]{$\beta_3$}}%
\put(59.3000,-28.0500){\makebox(0,0)[lb]{$\beta_4$}}%
\put(37.4000,-32.6400){\makebox(0,0)[lb]{$\gamma_1:=\alpha_4-\beta_4$}}%
%
\special{pn 20}%
\special{ar 2150 812 240 175  1.5707963 6.2831853}%
\special{ar 2150 812 240 175  0.0000000 1.5291537}%
%
\special{pn 20}%
\special{ar 2150 1643 240 175  1.5707963 6.2831853}%
\special{ar 2150 1643 240 175  0.0000000 1.5291537}%
%
\special{pn 20}%
\special{ar 2940 1242 240 175  1.5707963 6.2831853}%
\special{ar 2940 1242 240 175  0.0000000 1.5291537}%
%
\special{pn 20}%
\special{pa 2360 885}%
\special{pa 2750 1118}%
\special{fp}%
%
\special{pn 20}%
\special{pa 2370 1577}%
\special{pa 2770 1366}%
\special{fp}%
\put(19.6000,-17.1000){\makebox(0,0)[lb]{$x-t$}}%
\put(28.2000,-13.0500){\makebox(0,0)[lb]{$y$}}%
\put(19.2000,-8.7300){\makebox(0,0)[lb]{$x-\infty$}}%
\put(19.8000,-14.5800){\makebox(0,0)[lb]{$\alpha_0$}}%
\put(28.0000,-10.1700){\makebox(0,0)[lb]{$\alpha_2$}}%
%
\special{pn 8}%
\special{pa 3130 1131}%
\special{pa 3650 888}%
\special{dt 0.045}%
\special{pa 3650 888}%
\special{pa 3649 888}%
\special{dt 0.045}%
\special{pa 3100 1365}%
\special{pa 3670 1581}%
\special{dt 0.045}%
\special{pa 3670 1581}%
\special{pa 3669 1581}%
\special{dt 0.045}%
%
\special{pn 8}%
\special{ar 3890 816 240 176  1.5707963 1.6284886}%
\special{ar 3890 816 240 176  1.8015656 1.8592579}%
\special{ar 3890 816 240 176  2.0323348 2.0900271}%
\special{ar 3890 816 240 176  2.2631040 2.3207963}%
\special{ar 3890 816 240 176  2.4938732 2.5515656}%
\special{ar 3890 816 240 176  2.7246425 2.7823348}%
\special{ar 3890 816 240 176  2.9554117 3.0131040}%
\special{ar 3890 816 240 176  3.1861809 3.2438732}%
\special{ar 3890 816 240 176  3.4169502 3.4746425}%
\special{ar 3890 816 240 176  3.6477194 3.7054117}%
\special{ar 3890 816 240 176  3.8784886 3.9361809}%
\special{ar 3890 816 240 176  4.1092579 4.1669502}%
\special{ar 3890 816 240 176  4.3400271 4.3977194}%
\special{ar 3890 816 240 176  4.5707963 4.6284886}%
\special{ar 3890 816 240 176  4.8015656 4.8592579}%
\special{ar 3890 816 240 176  5.0323348 5.0900271}%
\special{ar 3890 816 240 176  5.2631040 5.3207963}%
\special{ar 3890 816 240 176  5.4938732 5.5515656}%
\special{ar 3890 816 240 176  5.7246425 5.7823348}%
\special{ar 3890 816 240 176  5.9554117 6.0131040}%
\special{ar 3890 816 240 176  6.1861809 6.2438732}%
\special{ar 3890 816 240 176  6.4169502 6.4746425}%
\special{ar 3890 816 240 176  6.6477194 6.7054117}%
\special{ar 3890 816 240 176  6.8784886 6.9361809}%
\special{ar 3890 816 240 176  7.1092579 7.1669502}%
\special{ar 3890 816 240 176  7.3400271 7.3977194}%
\special{ar 3890 816 240 176  7.5707963 7.6284886}%
\special{ar 3890 816 240 176  7.8015656 7.8123391}%
%
\special{pn 8}%
\special{ar 3880 1662 240 176  1.5707963 1.6284886}%
\special{ar 3880 1662 240 176  1.8015656 1.8592579}%
\special{ar 3880 1662 240 176  2.0323348 2.0900271}%
\special{ar 3880 1662 240 176  2.2631040 2.3207963}%
\special{ar 3880 1662 240 176  2.4938732 2.5515656}%
\special{ar 3880 1662 240 176  2.7246425 2.7823348}%
\special{ar 3880 1662 240 176  2.9554117 3.0131040}%
\special{ar 3880 1662 240 176  3.1861809 3.2438732}%
\special{ar 3880 1662 240 176  3.4169502 3.4746425}%
\special{ar 3880 1662 240 176  3.6477194 3.7054117}%
\special{ar 3880 1662 240 176  3.8784886 3.9361809}%
\special{ar 3880 1662 240 176  4.1092579 4.1669502}%
\special{ar 3880 1662 240 176  4.3400271 4.3977194}%
\special{ar 3880 1662 240 176  4.5707963 4.6284886}%
\special{ar 3880 1662 240 176  4.8015656 4.8592579}%
\special{ar 3880 1662 240 176  5.0323348 5.0900271}%
\special{ar 3880 1662 240 176  5.2631040 5.3207963}%
\special{ar 3880 1662 240 176  5.4938732 5.5515656}%
\special{ar 3880 1662 240 176  5.7246425 5.7823348}%
\special{ar 3880 1662 240 176  5.9554117 6.0131040}%
\special{ar 3880 1662 240 176  6.1861809 6.2438732}%
\special{ar 3880 1662 240 176  6.4169502 6.4746425}%
\special{ar 3880 1662 240 176  6.6477194 6.7054117}%
\special{ar 3880 1662 240 176  6.8784886 6.9361809}%
\special{ar 3880 1662 240 176  7.1092579 7.1669502}%
\special{ar 3880 1662 240 176  7.3400271 7.3977194}%
\special{ar 3880 1662 240 176  7.5707963 7.6284886}%
\special{ar 3880 1662 240 176  7.8015656 7.8123391}%
\put(37.1000,-9.0600){\makebox(0,0)[lb]{$x$}}%
\put(37.0000,-17.4300){\makebox(0,0)[lb]{$x-1$}}%
%
\special{pn 20}%
\special{ar 5690 1302 240 176  1.5707963 6.2831853}%
\special{ar 5690 1302 240 176  0.0000000 1.5291537}%
%
\special{pn 20}%
\special{ar 6580 865 240 175  1.5707963 6.2831853}%
\special{ar 6580 865 240 175  0.0000000 1.5291537}%
%
\special{pn 20}%
\special{ar 6580 1703 240 175  1.5707963 6.2831853}%
\special{ar 6580 1703 240 175  0.0000000 1.5291537}%
%
\special{pn 20}%
\special{pa 5900 1193}%
\special{pa 6370 974}%
\special{fp}%
%
\special{pn 20}%
\special{pa 5870 1412}%
\special{pa 6350 1645}%
\special{fp}%
\put(55.7000,-14.0100){\makebox(0,0)[lb]{$w$}}%
\put(63.8000,-17.7900){\makebox(0,0)[lb]{$z-1$}}%
\put(64.5000,-9.6900){\makebox(0,0)[lb]{$z$}}%
\put(55.1000,-10.9500){\makebox(0,0)[lb]{$\beta_2$}}%
\put(63.9000,-15.0900){\makebox(0,0)[lb]{$\beta_3$}}%
\put(64.1000,-6.6300){\makebox(0,0)[lb]{$\beta_4$}}%
%
\special{pn 8}%
\special{pa 5470 1203}%
\special{pa 4990 996}%
\special{dt 0.045}%
\special{pa 4990 996}%
\special{pa 4991 996}%
\special{dt 0.045}%
%
\special{pn 8}%
\special{pa 5530 1428}%
\special{pa 5020 1635}%
\special{dt 0.045}%
\special{pa 5020 1635}%
\special{pa 5021 1635}%
\special{dt 0.045}%
%
\special{pn 8}%
\special{ar 4810 879 240 175  1.5707963 1.6286277}%
\special{ar 4810 879 240 175  1.8021216 1.8599530}%
\special{ar 4810 879 240 175  2.0334469 2.0912783}%
\special{ar 4810 879 240 175  2.2647722 2.3226036}%
\special{ar 4810 879 240 175  2.4960975 2.5539289}%
\special{ar 4810 879 240 175  2.7274228 2.7852542}%
\special{ar 4810 879 240 175  2.9587481 3.0165795}%
\special{ar 4810 879 240 175  3.1900734 3.2479048}%
\special{ar 4810 879 240 175  3.4213987 3.4792301}%
\special{ar 4810 879 240 175  3.6527240 3.7105554}%
\special{ar 4810 879 240 175  3.8840493 3.9418807}%
\special{ar 4810 879 240 175  4.1153746 4.1732060}%
\special{ar 4810 879 240 175  4.3466999 4.4045313}%
\special{ar 4810 879 240 175  4.5780252 4.6358566}%
\special{ar 4810 879 240 175  4.8093505 4.8671819}%
\special{ar 4810 879 240 175  5.0406758 5.0985072}%
\special{ar 4810 879 240 175  5.2720011 5.3298325}%
\special{ar 4810 879 240 175  5.5033264 5.5611578}%
\special{ar 4810 879 240 175  5.7346517 5.7924831}%
\special{ar 4810 879 240 175  5.9659770 6.0238084}%
\special{ar 4810 879 240 175  6.1973024 6.2551337}%
\special{ar 4810 879 240 175  6.4286277 6.4864590}%
\special{ar 4810 879 240 175  6.6599530 6.7177843}%
\special{ar 4810 879 240 175  6.8912783 6.9491096}%
\special{ar 4810 879 240 175  7.1226036 7.1804349}%
\special{ar 4810 879 240 175  7.3539289 7.4117602}%
\special{ar 4810 879 240 175  7.5852542 7.6430855}%
%
\special{pn 8}%
\special{ar 4800 1662 240 176  1.5707963 1.6284886}%
\special{ar 4800 1662 240 176  1.8015656 1.8592579}%
\special{ar 4800 1662 240 176  2.0323348 2.0900271}%
\special{ar 4800 1662 240 176  2.2631040 2.3207963}%
\special{ar 4800 1662 240 176  2.4938732 2.5515656}%
\special{ar 4800 1662 240 176  2.7246425 2.7823348}%
\special{ar 4800 1662 240 176  2.9554117 3.0131040}%
\special{ar 4800 1662 240 176  3.1861809 3.2438732}%
\special{ar 4800 1662 240 176  3.4169502 3.4746425}%
\special{ar 4800 1662 240 176  3.6477194 3.7054117}%
\special{ar 4800 1662 240 176  3.8784886 3.9361809}%
\special{ar 4800 1662 240 176  4.1092579 4.1669502}%
\special{ar 4800 1662 240 176  4.3400271 4.3977194}%
\special{ar 4800 1662 240 176  4.5707963 4.6284886}%
\special{ar 4800 1662 240 176  4.8015656 4.8592579}%
\special{ar 4800 1662 240 176  5.0323348 5.0900271}%
\special{ar 4800 1662 240 176  5.2631040 5.3207963}%
\special{ar 4800 1662 240 176  5.4938732 5.5515656}%
\special{ar 4800 1662 240 176  5.7246425 5.7823348}%
\special{ar 4800 1662 240 176  5.9554117 6.0131040}%
\special{ar 4800 1662 240 176  6.1861809 6.2438732}%
\special{ar 4800 1662 240 176  6.4169502 6.4746425}%
\special{ar 4800 1662 240 176  6.6477194 6.7054117}%
\special{ar 4800 1662 240 176  6.8784886 6.9361809}%
\special{ar 4800 1662 240 176  7.1092579 7.1669502}%
\special{ar 4800 1662 240 176  7.3400271 7.3977194}%
\special{ar 4800 1662 240 176  7.5707963 7.6284886}%
\special{ar 4800 1662 240 176  7.8015656 7.8123391}%
%
\special{pn 20}%
\special{pa 3710 1986}%
\special{pa 3710 2364}%
\special{fp}%
\special{sh 1}%
\special{pa 3710 2364}%
\special{pa 3730 2297}%
\special{pa 3710 2311}%
\special{pa 3690 2297}%
\special{pa 3710 2364}%
\special{fp}%
%
\special{pn 20}%
\special{pa 5200 1977}%
\special{pa 5200 2355}%
\special{fp}%
\special{sh 1}%
\special{pa 5200 2355}%
\special{pa 5220 2288}%
\special{pa 5200 2302}%
\special{pa 5180 2288}%
\special{pa 5200 2355}%
\special{fp}%
%
\put(48.2000,-8.7900){\makebox(0,0)[lb]{}}%
\put(46.0000,-9.1500){\makebox(0,0)[lb]{$z-\infty$}}%
\put(45.9000,-17.2500){\makebox(0,0)[lb]{$z-t$}}%
\put(19.9000,-6.2000){\makebox(0,0)[lb]{$\alpha_1$}}%
\end{picture}%
\label{fig:CPVI2}
\caption{Dynkin diagram of type $D_6^{(1)}$}
\end{figure}
Here, let us explain Figure 1:
\begin{enumerate}
\item The above pictures of Figure 1 denote the Dynkin diagram of type $D_4^{(1)}$. The symbol in each circle denotes the invariant divisor of the sixth Painlev\'e system.
\item The below picture of Figure 1 denotes the Dynkin diagram of type $D_6^{(1)}$. The symbol in each circle denotes the invariant divisor of the system \eqref{1}.
\end{enumerate}

Here, let us consider the following problem.

\begin{problem}\label{prob2}
Can we find an algebraic ordinary differential system with the following assumptions {\rm (A1),(A2)? \rm}
\end{problem}

\begin{assumption}
$(A1)$ The system is a polynomial Hamiltonian system with the Hamiltonian $H \in {\Bbb C}(t)[x,y,z,w,\alpha_0,\alpha_1,\dots,\alpha_6]$.

$(A2)$ The system is invariant under the transformations $s_0,s_1,\dots,s_6$ given in Theorem \ref{1.1}.
\end{assumption}

In the brute force approach to Problem \ref{prob2}, we must deal with polynomials $H$ in
$$
\alpha_0,\alpha_1,\ldots,\alpha_6,t,x,y,z,w.
$$
This approach, however, soon came to a deadlock because of technical difficulties. In this paper, we present a new approach from the viewpoint of holomorphy.
\begin{theorem}\label{1.2}
Let us consider a polynomial Hamiltonian system with Hamiltonian $H \in {\Bbb C}[x,y,z,w]$. We assume that

$(A1)$ $deg(H)=5$ with respect to $x,y,z,w$.

$(A2)$ This system becomes again a polynomial Hamiltonian system in each coordinate system $(x_i,y_i,z_i,w_i) \ (i=0,1,\dots ,6)${\rm : \rm}
\begin{align*}
&r_0:x_0=-((x-t)y-\alpha_0)y, \ y_0=1/y, \ z_0=z, \ w_0=w, \\
&r_1:x_1=1/x, \ y_1=-x(xy+\alpha_1+\alpha_2), \ z_1=z, \ w_1=w, \\
&r_2:x_2=1/x, \ y_2=-x(xy+\alpha_2), \ z_2=z, \ w_2=w, \\
&r_3:x_3=-((x-z)y-\alpha_3)y, \ y_3=1/y, \ z_3=z, \ w_3=y+w, \\
&r_4:x_4=x, \ y_4=y, \ z_4=1/z, \ w_4=-z(zw+\alpha_4), \\
&r_5:x_5=x, \ y_5=y, \ z_5=-((z-1)w-\alpha_5)w, \ w_5=1/w, \\
&r_6:x_6=x, \ y_6=y, \ z_6=-w(zw-\alpha_6), \ w_6=1/w.
\end{align*}
Then such a system coincides with the system \eqref{1} with the polynomial Hamiltonian
\begin{align*}
H &=H_{VI}(x,y,t;\alpha_0,\alpha_1,\alpha_2,\alpha_3+2\alpha_4+\alpha_5,\alpha_3+\alpha_6)\\
&+H_{VI}(z,w,t;\alpha_0+\alpha_3,\alpha_1+2\alpha_2+\alpha_3,\alpha_4,\alpha_5,\alpha_6)+\frac{2(x-t)yz\{(z-1)w+\alpha_4\}}{t(t-1)}.
\end{align*}
\end{theorem}
By this theorem, we can also recover the parameter's relation \eqref{3}.

We note that the condition $(A2)$ should be read that
\begin{align*}
&r_j(H) \quad (j=1,2,\ldots,6), \quad r_0(H-y)
\end{align*}
are polynomials with respect to $x_i,y_i,z_i,w_i$.

In this method, $\alpha_0,\alpha_1,\ldots,\alpha_6$ can be treated as parameters rather than variables. In the holomorphy requirement, we only need to consider polynomials in $x,y,z,w$. Hence, the number of unknown coefficients can be drastically reduced.

Theorems \ref{1.1}, \ref{1.2} and Proposition \ref{P1.1} can be cheched by a direct calculation, respectively.

In addition to Theorems \ref{1.1} and \ref{1.2}, we give an explicit description of a confluence process to the system of type $A_5^{(1)}$.
\begin{theorem}\label{1.3}
For the system \eqref{1} of type $D_6^{(1)}$, we make the change of parameters and variables
\begin{gather}
\begin{gathered}\label{11}
\alpha_0={\varepsilon}^{-1}, \quad \alpha_1=A_0, \quad \alpha_2=A_1, \quad \alpha_4-\beta_4=A_2,\\
\beta_2=A_3, \quad \beta_3=-{\varepsilon}^{-1}-(A_1+A_2+A_3-A_5), \quad \beta_4=A_4,
\end{gathered}\\
\begin{gathered}\label{12}
t=1-{\varepsilon}T, \quad x=\frac{X}{X-T}, \quad z=\frac{Z}{Z-T},\\
y=-\frac{(X-T)\{(X-T)Y+A_1\}}{T}, \quad w=-\frac{(Z-T)\{(Z-T)W+A_3\}}{T}
\end{gathered}
\end{gather}
from $\alpha_0,\alpha_1,\alpha_2,\gamma_1,\beta_2,\beta_3,\beta_4,t,x,y,z,w$ to $A_0,\dots ,A_5,\varepsilon,T,X,Y,Z,W$. Then the system \eqref{1} can also be written in the new variables $T,X,Y,Z,W$ and parameters $A_0,A_1,\dots ,A_5,\varepsilon$ as a Hamiltonian system. This new system tends to the system of type $A_5^{(1)}$ as $\varepsilon \rightarrow 0$.
\end{theorem}
Here, the system of type $A_5^{(1)}$ is explicitly given in Section 6.

By proving the following theorem, we see how the transformation in Theorem \ref{1.3} works on the B{\"a}cklund transformation group $W(D_6^{(1)})=<s_0,s_1,\dots ,s_6>$ described in Theorem \ref{1.1}.

\begin{theorem}\label{1.4}
For the transformations \eqref{11}, \eqref{12} given in Theorem \ref{1.3} we can choose a subgroup $W_{D_6^{(1)} \rightarrow A_5^{(1)}}$ of the B{\"a}cklund transformation group $W(D_6^{(1)})$ so that $W_{D_6^{(1)} \rightarrow A_5^{(1)}}$ converges to the B{\"a}cklund transformation group $W(A_5^{(1)})$ of the system \eqref{29} with the polynomial Hamiltonian \eqref{30} \rm{(see Section 6) \rm}.
\end{theorem}

\section{The systems of type $B_6^{(1)}$}
In this section, we find two types of a 6-parameter family of coupled Painlev\'e VI systems in dimension four with affine Weyl group symmetry of type $B_6^{(1)}$. Each of them is equivalent to a polynomial Hamiltonian system, however, each has a different representaion of type $B_6^{(1)}$. We also show that each of them is equivalent to the system \eqref{1} by a birational and symplectic transformation.

The first member is given by
\begin{equation}\label{13}
\frac{dx}{dt}=\frac{\partial H}{\partial y}, \ \ \frac{dy}{dt}=-\frac{\partial H}{\partial x}, \ \ \frac{dz}{dt}=\frac{\partial H}{\partial w}, \ \ \frac{dw}{dt}=-\frac{\partial H}{\partial z}
\end{equation}
with the polynomial Hamiltonian
\begin{equation}\label{14}
\begin{aligned}
&H ={\tilde{H}}_{VI}(x,y,t;2\alpha_0+\alpha_1,\alpha_1,\alpha_2,\alpha_3+2\alpha_4+\alpha_5,\alpha_3+\alpha_6)\\
&+H_{VI}(z,w,t;2\alpha_0+\alpha_1+2\alpha_2+\alpha_3,\alpha_1+\alpha_3,\alpha_4,\alpha_5,\alpha_6)\\
  & \hspace{2cm} +\frac{2xz\{(tx-1)y+t\alpha_2\}\{(z-1)w+\alpha_4\}}{t(t-1)}.\\
\end{aligned}
\end{equation}
Here $x,y,z$ and $w$ denote unknown complex variables, and $\alpha_0,\alpha_1,\dots ,\alpha_6$ are complex parameters satisfying the relation
\begin{equation}\label{15}
\begin{aligned}
2\alpha_0+2\alpha_1+2\alpha_2+2\alpha_3+2\alpha_4+\alpha_5+\alpha_6=1.
\end{aligned}
\end{equation}

The symbol ${\tilde{H}}_{VI}(Q,P,t;\gamma_0,\gamma_1,\gamma_2,\gamma_3,\gamma_4)$ denotes the Hamiltonian of the second-order Painlev\'e VI equation given by
\begin{align}\label{16}
\begin{split}
&{\tilde{H}}_{VI}(Q,P,t;\gamma_0,\gamma_1,\gamma_2,\gamma_3,\gamma_4)\\
&=\frac{1}{t(t-1)}[P^2(tQ-1)(Q-1)Q-\{(\gamma_0-1)t(Q-1)Q+\gamma_1(Q-1)(tQ-1)\\
&+\gamma_3Q(tQ-1)\}P+\gamma_2(\gamma_2+\gamma_4)tQ]  \ \ (\gamma_0+\gamma_1+2\gamma_2+\gamma_3+\gamma_4=1).
\end{split}
\end{align}

\begin{theorem}\label{2.1}
The system \eqref{13} admits extended affine Weyl group symmetry of type $B_6^{(1)}$ as the group of its B{\"a}cklund transformations, whose generators $S_0,S_1,\dots ,S_6,\varphi$ are explicitly given as follows{\rm : \rm}with the notation$:(*):=(x,y,z,w,t;\alpha_0,\alpha_1,\ldots,\alpha_6),$
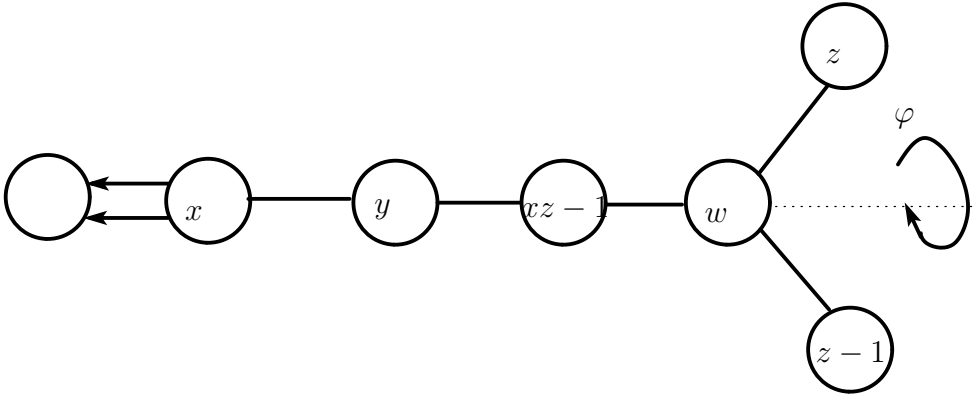
\begin{figure}[h]
\unitlength 0.1in
\begin{picture}(51.00,20.30)(20.50,-29.80)
%
\special{pn 20}%
\special{ar 2270 1960 220 220  1.6078164 6.2831853}%
\special{ar 2270 1960 220 220  0.0000000 1.6064954}%
%
\special{pn 20}%
\special{ar 3110 1980 220 220  1.6078164 6.2831853}%
\special{ar 3110 1980 220 220  0.0000000 1.6064954}%
%
\special{pn 20}%
\special{ar 4090 1990 220 220  1.6078164 6.2831853}%
\special{ar 4090 1990 220 220  0.0000000 1.6064954}%
%
\special{pn 20}%
\special{ar 4970 1990 220 220  1.6078164 6.2831853}%
\special{ar 4970 1990 220 220  0.0000000 1.6064954}%
%
\special{pn 20}%
\special{ar 5830 1990 220 220  1.6078164 6.2831853}%
\special{ar 5830 1990 220 220  0.0000000 1.6064954}%
%
\special{pn 20}%
\special{ar 6440 1170 220 220  1.6078164 6.2831853}%
\special{ar 6440 1170 220 220  0.0000000 1.6064954}%
%
\special{pn 20}%
\special{ar 6470 2760 220 220  1.6078164 6.2831853}%
\special{ar 6470 2760 220 220  0.0000000 1.6064954}%
%
\special{pn 20}%
\special{pa 3320 1970}%
\special{pa 3850 1970}%
\special{fp}%
%
\special{pn 20}%
\special{pa 4320 1990}%
\special{pa 4740 1990}%
\special{fp}%
%
\special{pn 20}%
\special{pa 5190 2000}%
\special{pa 5590 2000}%
\special{fp}%
%
\special{pn 20}%
\special{pa 5990 1840}%
\special{pa 6350 1380}%
\special{fp}%
%
\special{pn 20}%
\special{pa 6000 2130}%
\special{pa 6360 2550}%
\special{fp}%
%
\special{pn 20}%
\special{pa 2900 1890}%
\special{pa 2500 1890}%
\special{fp}%
\special{sh 1}%
\special{pa 2500 1890}%
\special{pa 2567 1910}%
\special{pa 2553 1890}%
\special{pa 2567 1870}%
\special{pa 2500 1890}%
\special{fp}%
%
\special{pn 20}%
\special{pa 2900 2070}%
\special{pa 2490 2070}%
\special{fp}%
\special{sh 1}%
\special{pa 2490 2070}%
\special{pa 2557 2090}%
\special{pa 2543 2070}%
\special{pa 2557 2050}%
\special{pa 2490 2070}%
\special{fp}%
\put(29.9000,-20.8000){\makebox(0,0)[lb]{$x$}}%
\put(39.8000,-20.8000){\makebox(0,0)[lb]{$y$}}%
\put(47.5000,-20.8000){\makebox(0,0)[lb]{$xz-1$}}%
\put(57.1000,-20.9000){\makebox(0,0)[lb]{$w$}}%
\put(63.4000,-12.6000){\makebox(0,0)[lb]{$z$}}%
\put(62.9000,-28.4000){\makebox(0,0)[lb]{$z-1$}}%
%
\special{pn 8}%
\special{pa 6080 2010}%
\special{pa 7150 2010}%
\special{dt 0.045}%
\special{pa 7150 2010}%
\special{pa 7149 2010}%
\special{dt 0.045}%
%
\special{pn 20}%
\special{pa 6720 1790}%
\special{pa 6741 1757}%
\special{pa 6763 1725}%
\special{pa 6785 1697}%
\special{pa 6808 1674}%
\special{pa 6831 1658}%
\special{pa 6855 1650}%
\special{pa 6881 1653}%
\special{pa 6907 1665}%
\special{pa 6933 1684}%
\special{pa 6959 1708}%
\special{pa 6983 1737}%
\special{pa 7005 1768}%
\special{pa 7025 1799}%
\special{pa 7042 1831}%
\special{pa 7057 1864}%
\special{pa 7068 1897}%
\special{pa 7077 1929}%
\special{pa 7083 1962}%
\special{pa 7087 1995}%
\special{pa 7087 2027}%
\special{pa 7084 2059}%
\special{pa 7077 2090}%
\special{pa 7068 2120}%
\special{pa 7055 2150}%
\special{pa 7039 2178}%
\special{pa 7018 2203}%
\special{pa 6990 2220}%
\special{pa 6954 2227}%
\special{pa 6916 2223}%
\special{pa 6880 2210}%
\special{pa 6853 2189}%
\special{pa 6840 2162}%
\special{pa 6840 2150}%
\special{sp}%
%
\special{pn 20}%
\special{pa 6840 2170}%
\special{pa 6770 2020}%
\special{fp}%
\special{sh 1}%
\special{pa 6770 2020}%
\special{pa 6780 2089}%
\special{pa 6793 2068}%
\special{pa 6816 2072}%
\special{pa 6770 2020}%
\special{fp}%
\put(67.0000,-16.0000){\makebox(0,0)[lb]{$\varphi$}}%
\end{picture}%
\label{fig:CPVI20}
\caption{Dynkin diagram of type $B_6^{(1)}$}
\end{figure}
\begin{align*}
\begin{split}
S_0: (*) &\rightarrow (\frac{tx-1}{t-1},\frac{(t-1)y}{t},\frac{(t-1)z}{t-z},\frac{(t-z)(tw-zw-\alpha_4)}{t(t-1)},1-t;\\
&-\alpha_0,\alpha_1+2\alpha_0,\alpha_2,\alpha_3,\alpha_4,\alpha_5,\alpha_6),
\end{split}\\
S_1: (*) &\rightarrow \left(x,y-\frac{\alpha_1}{x},z,w,t;\alpha_0+\alpha_1,-\alpha_1,\alpha_2+\alpha_1,\alpha_3,\alpha_4,\alpha_5,\alpha_6 \right), \\
S_2: (*) &\rightarrow \left(x+\frac{\alpha_2}{y},y,z,w,t;\alpha_0,\alpha_1+\alpha_2,-\alpha_2,\alpha_3+\alpha_2,\alpha_4,\alpha_5,\alpha_6 \right), \\
S_3: (*) &\rightarrow \left(x,y-\frac{\alpha_3z}{xz-1},z,w-\frac{\alpha_3x}{xz-1},t;\alpha_0,\alpha_1,\alpha_2+\alpha_3,-\alpha_3,\alpha_4+\alpha_3,\alpha_5,\alpha_6 \right), \\
S_4: (*) &\rightarrow \left(x,y,z+\frac{\alpha_4}{w},w,t;\alpha_0,\alpha_1,\alpha_2,\alpha_3+\alpha_4,-\alpha_4,\alpha_5+\alpha_4,\alpha_6+\alpha_4 \right), \\
S_5: (*) &\rightarrow \left(x,y,z,w-\frac{\alpha_5}{z-1},t;\alpha_0,\alpha_1,\alpha_2,\alpha_3,\alpha_4+\alpha_5,-\alpha_5,\alpha_6 \right), \\
S_6: (*) &\rightarrow \left(x,y,z,w-\frac{\alpha_6}{z},t;\alpha_0,\alpha_1,\alpha_2,\alpha_3,\alpha_4+\alpha_6,\alpha_5,-\alpha_6 \right), \\
\begin{split}
\varphi: (*) &\rightarrow  (\frac{x}{x-1},-(x-1)\{(x-1)y+\alpha_2\},1-z,-w,1-t;\\
& \qquad \alpha_0,\alpha_1,\alpha_2,\alpha_3,\alpha_4,\alpha_6,\alpha_5).
\end{split}
\end{align*}
\end{theorem}
We remark that the B{\"a}cklund transformations $S_1,\ldots,S_6$ satisfy the relation \eqref{unvre}. However, the transformation $S_0$ do not satisfy so.

Theorem \ref{2.1} can be cheched by a direct calculation.

\begin{theorem}\label{2.2}
For the system \eqref{1} of type $D_6^{(1)}$, we make the change of parameters and variables
\begin{gather}
\begin{gathered}\label{17}
\frac{(\alpha_0-\alpha_1)}{2}=A_0, \quad \alpha_1=A_1, \quad \alpha_2=A_2, \quad \alpha_4-\beta_4=A_3,\\
 \ \beta_2=A_4, \quad \beta_3=A_5, \quad \beta_4=A_6,
\end{gathered}\\
\begin{gathered}\label{18}
X=\frac{1}{x}, \quad Y=-(xy+\alpha_2)x, \quad Z=z, \quad W=w
\end{gathered}
\end{gather}
from $\alpha_0,\alpha_1,\alpha_2,\alpha_4,\beta_2,\beta_3,\beta_4,x,y,z,w$ to $A_0,A_1,\dots ,A_6,X,Y,Z,W$. Then the system \eqref{1} can also be written in the new variables $X,Y,Z,W$ and parameters $A_0,\dots ,A_6$ as a Hamiltonian system. This new system tends to the system \eqref{13} with the Hamiltonian \eqref{14}.
\end{theorem}

\begin{proof}
Notice that
$$
2A_0+2A_1+2A_2+2A_3+2A_4+A_5+A_6=\alpha_0+\alpha_1+2\alpha_2+2(\alpha_4-\beta_4)+2\beta_2+\beta_3+\beta_4=1
$$
and the change of variables from $(x,y,z,w)$ to $(X,Y,Z,W)$ is symplectic. Choose $S_i \ (i=0,1,\dots ,6)$ and $\varphi$ as
$$
S_0:=\pi_4, \ S_1:=s_1, \ S_2:=s_2, \ S_3:=s_3, \ S_4:=s_4, \ S_5:=s_5, \ S_6:=s_6, \ \varphi:=\pi_3.
$$
Then the transformations $S_i$ are reflections of the parameters $A_0,A_1,\dots ,A_6$. The transformation group $\tilde{W}(B_6^{(1)})=<S_0,S_1,\dots ,S_6,\varphi>$ coincides with the transformations given in Theorem \ref{2.1}.
\end{proof}

The second member is given by
\begin{equation}\label{19}
\frac{dx}{dt}=\frac{\partial H}{\partial y}, \ \ \frac{dy}{dt}=-\frac{\partial H}{\partial x}, \ \ \frac{dz}{dt}=\frac{\partial H}{\partial w}, \ \ \frac{dw}{dt}=-\frac{\partial H}{\partial z}
\end{equation}
with the polynomial Hamiltonian
\begin{equation}\label{20}
\begin{aligned}
H &=H_{VI}(x,y,t;\alpha_0,\alpha_1,\alpha_2,\alpha_3+\alpha_5+2\alpha_6,\alpha_3+2\alpha_4+\alpha_5)\\
&+\tilde{H}_{VI}(z,w,t;\alpha_0+\alpha_3,\alpha_1+2\alpha_2+\alpha_3,\alpha_4,\alpha_5+2\alpha_6,\alpha_5)\\
  & \hspace{2cm} +\frac{2(x-t)y(z-1)w}{t(t-1)}.\\
\end{aligned}
\end{equation}
Here $x,y,z$ and $w$ denote unknown complex variables, and $\alpha_0,\alpha_1,\dots ,\alpha_6$ are complex parameters satisfying the relation
\begin{equation}
\begin{aligned}
\alpha_0+\alpha_1+2\alpha_2+2\alpha_3+2\alpha_4+2\alpha_5+2\alpha_6=1.
\end{aligned}
\end{equation}

\begin{theorem}\label{2.3}
The system \eqref{19} admits extended affine Weyl group symmetry of type $B_6^{(1)}$ as the group of its B{\"a}cklund transformations, whose generators $w_0,w_1,w_2,\dots ,w_6,\psi$ are explicitly given as follows{\rm : \rm}with the notation$:(*):=(x,y,z,w,t;\alpha_0,\alpha_1,\ldots,\alpha_6),$
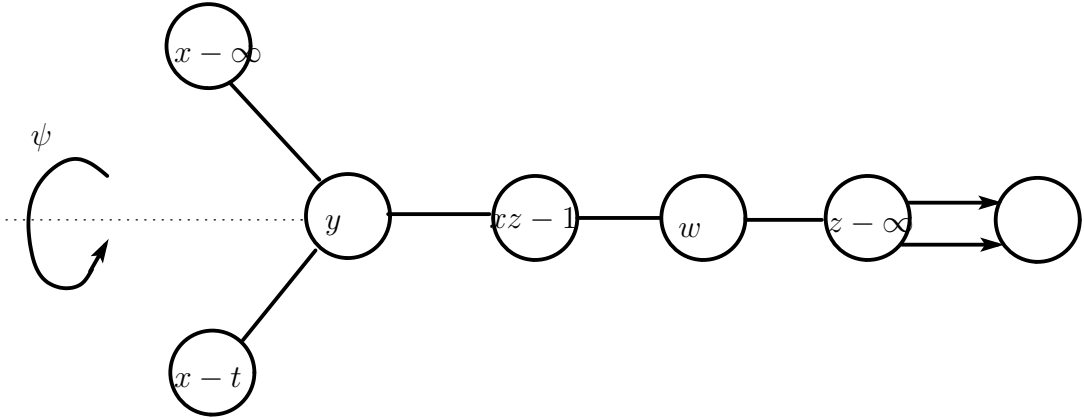
\begin{figure}[h]
\unitlength 0.1in
\begin{picture}(56.20,21.50)(8.60,-30.30)
%
\special{pn 20}%
\special{ar 2650 1990 220 220  1.6078164 6.2831853}%
\special{ar 2650 1990 220 220  0.0000000 1.6064954}%
%
\special{pn 20}%
\special{ar 3630 2000 220 220  1.6078164 6.2831853}%
\special{ar 3630 2000 220 220  0.0000000 1.6064954}%
%
\special{pn 20}%
\special{ar 4510 2000 220 220  1.6078164 6.2831853}%
\special{ar 4510 2000 220 220  0.0000000 1.6064954}%
%
\special{pn 20}%
\special{ar 5370 2000 220 220  1.6078164 6.2831853}%
\special{ar 5370 2000 220 220  0.0000000 1.6064954}%
%
\special{pn 20}%
\special{pa 2860 1980}%
\special{pa 3390 1980}%
\special{fp}%
%
\special{pn 20}%
\special{pa 3860 2000}%
\special{pa 4280 2000}%
\special{fp}%
%
\special{pn 20}%
\special{pa 4730 2010}%
\special{pa 5130 2010}%
\special{fp}%
\put(25.3000,-20.9000){\makebox(0,0)[lb]{$y$}}%
\put(33.9000,-20.7000){\makebox(0,0)[lb]{$xz-1$}}%
\put(43.8000,-21.0000){\makebox(0,0)[lb]{$w$}}%
\put(51.6000,-20.9000){\makebox(0,0)[lb]{$z-\infty$}}%
%
\special{pn 20}%
\special{ar 1920 1100 220 220  1.6078164 6.2831853}%
\special{ar 1920 1100 220 220  0.0000000 1.6064954}%
%
\special{pn 20}%
\special{ar 1940 2810 220 220  1.6078164 6.2831853}%
\special{ar 1940 2810 220 220  0.0000000 1.6064954}%
%
\special{pn 20}%
\special{ar 6260 2010 220 220  1.6078164 6.2831853}%
\special{ar 6260 2010 220 220  0.0000000 1.6064954}%
%
\special{pn 20}%
\special{pa 2100 2640}%
\special{pa 2480 2170}%
\special{fp}%
%
\special{pn 20}%
\special{pa 2030 1290}%
\special{pa 2500 1800}%
\special{fp}%
%
\special{pn 20}%
\special{pa 5580 1920}%
\special{pa 6040 1920}%
\special{fp}%
\special{sh 1}%
\special{pa 6040 1920}%
\special{pa 5973 1900}%
\special{pa 5987 1920}%
\special{pa 5973 1940}%
\special{pa 6040 1920}%
\special{fp}%
%
\special{pn 20}%
\special{pa 5550 2140}%
\special{pa 6050 2140}%
\special{fp}%
\special{sh 1}%
\special{pa 6050 2140}%
\special{pa 5983 2120}%
\special{pa 5997 2140}%
\special{pa 5983 2160}%
\special{pa 6050 2140}%
\special{fp}%
\put(17.4000,-29.0000){\makebox(0,0)[lb]{$x-t$}}%
\put(17.4000,-12.0000){\makebox(0,0)[lb]{$x-\infty$}}%
%
\special{pn 8}%
\special{pa 860 2010}%
\special{pa 2410 2010}%
\special{dt 0.045}%
\special{pa 2410 2010}%
\special{pa 2409 2010}%
\special{dt 0.045}%
%
\special{pn 20}%
\special{pa 1390 1780}%
\special{pa 1363 1757}%
\special{pa 1336 1736}%
\special{pa 1308 1717}%
\special{pa 1280 1703}%
\special{pa 1251 1693}%
\special{pa 1222 1690}%
\special{pa 1192 1694}%
\special{pa 1163 1704}%
\special{pa 1134 1720}%
\special{pa 1106 1740}%
\special{pa 1080 1763}%
\special{pa 1056 1789}%
\special{pa 1035 1817}%
\special{pa 1018 1845}%
\special{pa 1005 1875}%
\special{pa 994 1905}%
\special{pa 986 1936}%
\special{pa 981 1968}%
\special{pa 978 2000}%
\special{pa 977 2033}%
\special{pa 978 2066}%
\special{pa 981 2099}%
\special{pa 985 2133}%
\special{pa 991 2166}%
\special{pa 998 2200}%
\special{pa 1007 2233}%
\special{pa 1019 2264}%
\special{pa 1034 2292}%
\special{pa 1053 2317}%
\special{pa 1077 2338}%
\special{pa 1107 2354}%
\special{pa 1140 2365}%
\special{pa 1174 2369}%
\special{pa 1208 2366}%
\special{pa 1239 2356}%
\special{pa 1264 2339}%
\special{pa 1284 2314}%
\special{pa 1301 2286}%
\special{pa 1310 2270}%
\special{sp}%
%
\special{pn 20}%
\special{pa 1300 2280}%
\special{pa 1380 2140}%
\special{fp}%
\special{sh 1}%
\special{pa 1380 2140}%
\special{pa 1330 2188}%
\special{pa 1354 2186}%
\special{pa 1364 2208}%
\special{pa 1380 2140}%
\special{fp}%
\put(9.9000,-16.5000){\makebox(0,0)[lb]{$\psi$}}%
\end{picture}%
\label{fig:CPVI21}
\caption{Dynkin diagram of type $B_6^{(1)}$}
\end{figure}
\begin{align*}
w_0: (*) &\rightarrow \left(x,y-\frac{\alpha_0}{x-t},z,w,t;-\alpha_0,\alpha_1,\alpha_2+\alpha_0,\alpha_3,\alpha_4,\alpha_5,\alpha_6 \right), \\
w_1: (*) &\rightarrow (x,y,z,w,t;\alpha_0,-\alpha_1,\alpha_2+\alpha_1,\alpha_3,\alpha_4,\alpha_5,\alpha_6),\\
w_2: (*) &\rightarrow \left(x+\frac{\alpha_2}{y},y,z,w,t;\alpha_0+\alpha_2,\alpha_1+\alpha_2,-\alpha_2,\alpha_3+\alpha_2,\alpha_4,\alpha_5,\alpha_6 \right), \\
w_3: (*) &\rightarrow \left(x,y-\frac{\alpha_3z}{xz-1},z,w-\frac{\alpha_3x}{xz-1},t;\alpha_0,\alpha_1,\alpha_2+\alpha_3,-\alpha_3,\alpha_4+\alpha_3,\alpha_5,\alpha_6 \right), \\
w_4: (*) &\rightarrow \left(x,y,z+\frac{\alpha_4}{w},w,t;\alpha_0,\alpha_1,\alpha_2,\alpha_3+\alpha_4,-\alpha_4,\alpha_5+\alpha_4,\alpha_6 \right), \\
w_5: (*) &\rightarrow (x,y,z,w,t;\alpha_0,\alpha_1,\alpha_2,\alpha_3,\alpha_4+\alpha_5,-\alpha_5,\alpha_6+\alpha_5), \\
\begin{split}
w_6: (*) &\rightarrow (1-x,-y,\frac{z}{z-1},-(z-1)\{(z-1)w+\alpha_4\},1-t; \\
&\alpha_0,\alpha_1,\alpha_2,\alpha_3,\alpha_4,\alpha_5+2\alpha_6,-\alpha_6),
\end{split}\\
\begin{split}
\psi: (*) &\rightarrow  (\frac{(t-1)x}{t-x},\frac{(t-x)\{(t-x)y-\alpha_2\}}{t(t-1)},\frac{tz-1}{t-1},\frac{(t-1)w}{t},1-t;\\
&\alpha_1,\alpha_0,\alpha_2,\alpha_3,\alpha_4,\alpha_5,\alpha_6).
\end{split}
\end{align*}
\end{theorem}
We remark that the B{\"a}cklund transformations $w_0,\ldots,w_5$ satisfy the relation \eqref{unvre}. However, the transformation $w_6$ do not satisfy so.

Theorem \ref{2.3} can be cheched by a direct calculation.

\begin{theorem}\label{2.4}
For the system \eqref{1} of type $D_6^{(1)}$, we make the change of parameters and variables
\begin{gather}
\begin{gathered}\label{22}
\alpha_0=A_0, \quad \alpha_1=A_1, \quad \alpha_2=A_2, \quad \alpha_4-\beta_4=A_3,\\
 \quad \beta_2=A_4, \quad \beta_4=A_5, \quad \frac{(\beta_3-\beta_4)}{2}=A_6,
\end{gathered}\\
\begin{gathered}\label{23}
X=x, \quad Y=y, \quad Z=\frac{1}{z}, \quad W=-z(zw+\beta_2)
\end{gathered}
\end{gather}
from $\alpha_0,\alpha_1,\alpha_2,\alpha_4,\beta_2,\beta_3,\beta_4,x,y,z,w$ to $A_0,A_1,\dots ,A_6,X,Y,Z,W$. Then the system \eqref{1} can also be written in the new variables $X,Y,Z,W$ and parameters $A_0,A_1,\dots ,A_6$ as a Hamiltonian system. This new system tends to the system \eqref{19} with the Hamiltonian \eqref{20}.
\end{theorem}

{\it Proof.} \qquad Notice that 
$$
A_0+A_1+2A_2+2A_3+2A_4+2A_5+2A_6=\alpha_0+\alpha_1+2\alpha_2+2(\alpha_4-\beta_4)+2\beta_2+\beta_3+\beta_4=1
$$
and the change of variables from $(x,y,z,w)$ to $(X,Y,Z,W)$ is symplectic. Choose $w_i \ (i=0,1,\dots ,6)$ and $\psi$ as
$$
w_0:=s_0, \ w_1:=s_1, \ w_2:=s_2, \ w_3:=s_3, \ w_4:=s_4, \ w_5:=s_6, \ w_6:=\pi_3, \ \psi:=\pi_4.
$$
Then the transformations $w_i$ are reflections of the parameters $A_0,A_1,\dots ,A_6$. The transformation group $\tilde{W}(B_6^{(1)})=<w_0,w_1,\dots ,w_6,\psi>$ coincides with the transformations given in Theorem \ref{2.3}. \qed

\section{The system of type $D_7^{(2)}$}
In this section, we find a 6-parameter family of coupled Painlev\'e VI systems in dimension four with affine Weyl group symmetry of type $D_7^{(2)}$. This system is equivalent to a polynomial Hamiltonian system. In the final stage of this section, this system is equivalent to the system \eqref{1} by a birational and symplectic transformation.

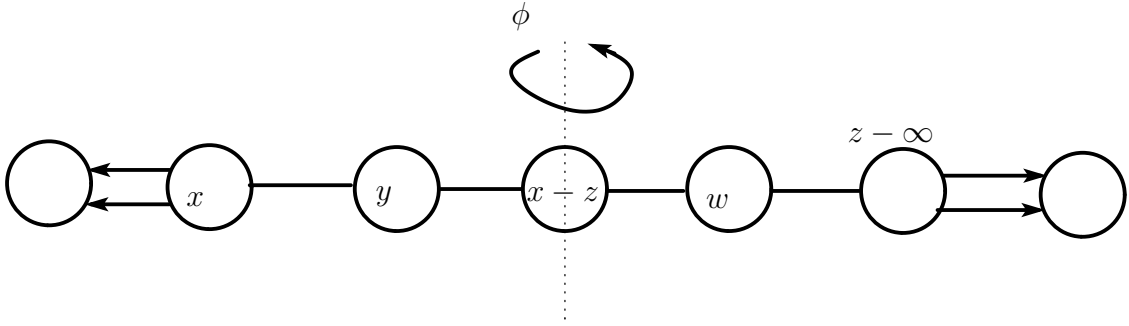
\begin{figure}
\unitlength 0.1in
\begin{picture}(58.50,16.80)(16.60,-27.90)
%
\special{pn 20}%
\special{ar 1880 2080 220 220  1.6078164 6.2831853}%
\special{ar 1880 2080 220 220  0.0000000 1.6064954}%
%
\special{pn 20}%
\special{ar 2720 2100 220 220  1.6078164 6.2831853}%
\special{ar 2720 2100 220 220  0.0000000 1.6064954}%
%
\special{pn 20}%
\special{ar 3700 2110 220 220  1.6078164 6.2831853}%
\special{ar 3700 2110 220 220  0.0000000 1.6064954}%
%
\special{pn 20}%
\special{ar 4580 2110 220 220  1.6078164 6.2831853}%
\special{ar 4580 2110 220 220  0.0000000 1.6064954}%
%
\special{pn 20}%
\special{ar 5440 2110 220 220  1.6078164 6.2831853}%
\special{ar 5440 2110 220 220  0.0000000 1.6064954}%
%
\special{pn 20}%
\special{pa 2930 2090}%
\special{pa 3460 2090}%
\special{fp}%
%
\special{pn 20}%
\special{pa 3930 2110}%
\special{pa 4350 2110}%
\special{fp}%
%
\special{pn 20}%
\special{pa 4800 2120}%
\special{pa 5200 2120}%
\special{fp}%
%
\special{pn 20}%
\special{pa 2510 2010}%
\special{pa 2110 2010}%
\special{fp}%
\special{sh 1}%
\special{pa 2110 2010}%
\special{pa 2177 2030}%
\special{pa 2163 2010}%
\special{pa 2177 1990}%
\special{pa 2110 2010}%
\special{fp}%
%
\special{pn 20}%
\special{pa 2510 2190}%
\special{pa 2100 2190}%
\special{fp}%
\special{sh 1}%
\special{pa 2100 2190}%
\special{pa 2167 2210}%
\special{pa 2153 2190}%
\special{pa 2167 2170}%
\special{pa 2100 2190}%
\special{fp}%
\put(26.0000,-22.0000){\makebox(0,0)[lb]{$x$}}%
\put(35.9000,-22.0000){\makebox(0,0)[lb]{$y$}}%
\put(53.2000,-22.1000){\makebox(0,0)[lb]{$w$}}%
%
\special{pn 20}%
\special{ar 6350 2120 220 220  1.6078164 6.2831853}%
\special{ar 6350 2120 220 220  0.0000000 1.6064954}%
%
\special{pn 20}%
\special{ar 7290 2140 220 220  1.6078164 6.2831853}%
\special{ar 7290 2140 220 220  0.0000000 1.6064954}%
%
\special{pn 20}%
\special{pa 6560 2040}%
\special{pa 7070 2040}%
\special{fp}%
\special{sh 1}%
\special{pa 7070 2040}%
\special{pa 7003 2020}%
\special{pa 7017 2040}%
\special{pa 7003 2060}%
\special{pa 7070 2040}%
\special{fp}%
%
\special{pn 20}%
\special{pa 6530 2220}%
\special{pa 7060 2220}%
\special{fp}%
\special{sh 1}%
\special{pa 7060 2220}%
\special{pa 6993 2200}%
\special{pa 7007 2220}%
\special{pa 6993 2240}%
\special{pa 7060 2220}%
\special{fp}%
%
\special{pn 20}%
\special{pa 5660 2120}%
\special{pa 6120 2120}%
\special{fp}%
\put(43.8000,-21.9000){\makebox(0,0)[lb]{$x-z$}}%
\put(60.6000,-18.8000){\makebox(0,0)[lb]{$z-\infty$}}%
%
\special{pn 8}%
\special{pa 4580 1340}%
\special{pa 4580 2790}%
\special{dt 0.045}%
\special{pa 4580 2790}%
\special{pa 4580 2789}%
\special{dt 0.045}%
%
\special{pn 20}%
\special{pa 4440 1390}%
\special{pa 4404 1406}%
\special{pa 4370 1422}%
\special{pa 4340 1439}%
\special{pa 4317 1457}%
\special{pa 4302 1476}%
\special{pa 4297 1497}%
\special{pa 4302 1519}%
\special{pa 4316 1541}%
\special{pa 4338 1564}%
\special{pa 4365 1587}%
\special{pa 4397 1609}%
\special{pa 4433 1630}%
\special{pa 4470 1649}%
\special{pa 4509 1666}%
\special{pa 4547 1680}%
\special{pa 4585 1691}%
\special{pa 4621 1699}%
\special{pa 4657 1704}%
\special{pa 4691 1705}%
\special{pa 4724 1703}%
\special{pa 4755 1698}%
\special{pa 4785 1688}%
\special{pa 4814 1675}%
\special{pa 4840 1658}%
\special{pa 4865 1636}%
\special{pa 4888 1610}%
\special{pa 4908 1581}%
\special{pa 4923 1550}%
\special{pa 4931 1519}%
\special{pa 4930 1490}%
\special{pa 4918 1465}%
\special{pa 4898 1443}%
\special{pa 4871 1423}%
\special{pa 4839 1405}%
\special{pa 4810 1390}%
\special{sp}%
%
\special{pn 20}%
\special{pa 4840 1400}%
\special{pa 4730 1360}%
\special{fp}%
\special{sh 1}%
\special{pa 4730 1360}%
\special{pa 4786 1402}%
\special{pa 4780 1378}%
\special{pa 4799 1364}%
\special{pa 4730 1360}%
\special{fp}%
\put(43.0000,-12.8000){\makebox(0,0)[lb]{$\phi$}}%
\end{picture}%
\label{fig:CPVI22}
\caption{Dynkin diagram of type $D_7^{(2)}$}
\end{figure}

This system is explicitly given as follows:
\begin{equation}\label{24}
\frac{dx}{dt}=\frac{\partial H}{\partial y}, \ \ \frac{dy}{dt}=-\frac{\partial H}{\partial x}, \ \ \frac{dz}{dt}=\frac{\partial H}{\partial w}, \ \ \frac{dw}{dt}=-\frac{\partial H}{\partial z}
\end{equation}
with the polynomial Hamiltonian
\begin{equation}\label{25}
\begin{aligned}
H &={\tilde{H}}_{VI}(x,y,t;2\alpha_0+\alpha_1,\alpha_1,\alpha_2,\alpha_3+\alpha_5+2\alpha_6,\alpha_3+2\alpha_4+\alpha_5)\\
&+\tilde{H}_{VI}(z,w,t;2\alpha_0+\alpha_1+2\alpha_2+\alpha_3,\alpha_1+\alpha_3,\alpha_4,\alpha_5+2\alpha_6,\alpha_5)\\
  & \hspace{2cm} +\frac{2x\{(tx-1)y+t\alpha_2\}(z-1)w}{t(t-1)}.\\
\end{aligned}
\end{equation}
Here $x,y,z$ and $w$ denote unknown complex variables, and $\alpha_0,\alpha_1,\dots ,\alpha_6$ are complex parameters satisfying the relation
\begin{equation}\label{26}
\begin{aligned}
\alpha_0+\alpha_1+\alpha_2+\alpha_3+\alpha_4+\alpha_5+\alpha_6=\frac{1}{2}.
\end{aligned}
\end{equation}

\begin{theorem}\label{3.1}
The system \eqref{24} admits extended affine Weyl group symmetry of type $D_7^{(2)}$ as the group of its B{\"a}cklund transformations, whose generators $u_0,u_1,\dots ,u_6,\phi$ are explicitly given as follows{\rm : \rm}with the notation$:(*):=(x,y,z,w,t;\alpha_0,\alpha_1,\ldots,\alpha_6),$
\begin{align*}
\begin{split}
        u_0: (*) &\rightarrow (\frac{tx-1}{t-1},\frac{(t-1)y}{t},\frac{tz-1}{t-1},\frac{(t-1)w}{t},1-t;\\
        &\qquad -\alpha_0,\alpha_1+2\alpha_0,\alpha_2,\alpha_3,\alpha_4,\alpha_5,\alpha_6),
\end{split}\\
        u_1: (*) &\rightarrow \left(x,y-\frac{\alpha_1}{x},z,w,t;\alpha_0+\alpha_1,-\alpha_1,\alpha_2+\alpha_1,\alpha_3,\alpha_4,\alpha_5,\alpha_6 \right), \\
        u_2: (*) &\rightarrow \left(x+\frac{\alpha_2}{y},y,z,w,t;\alpha_0,\alpha_1+\alpha_2,-\alpha_2,\alpha_3+\alpha_2,\alpha_4,\alpha_5,\alpha_6 \right), \\
        u_3: (*) &\rightarrow \left(x,y-\frac{\alpha_3}{x-z},z,w+\frac{\alpha_3}{x-z},t;\alpha_0,\alpha_1,\alpha_2+\alpha_3,-\alpha_3,\alpha_4+\alpha_3,\alpha_5,\alpha_6 \right), \\
        u_4: (*) &\rightarrow \left(x,y,z+\frac{\alpha_4}{w},w,t;\alpha_0,\alpha_1,\alpha_2,\alpha_3+\alpha_4,-\alpha_4,\alpha_5+\alpha_4,\alpha_6 \right), \\
       u_5: (*) &\rightarrow (x,y,z,w,t;\alpha_0,\alpha_1,\alpha_2,\alpha_3,\alpha_4+\alpha_5,-\alpha_5,\alpha_6+\alpha_5), \\
\begin{split}
       u_6: (*) &\rightarrow (\frac{x}{x-1},-(x-1)\{(x-1)y+\alpha_2\},\frac{z}{z-1},-(z-1)\{(z-1)w+\alpha_4\},1-t;\\
       & \qquad \alpha_0,\alpha_1,\alpha_2,\alpha_3,\alpha_4,\alpha_5+2\alpha_6,-\alpha_6),
\end{split}\\
       \phi: (*) &\rightarrow  \left(\frac{1}{tz},-tz(zw+\alpha_4),\frac{1}{tx},-tx(xy+\alpha_2),t;\alpha_6,\alpha_5,\alpha_4,\alpha_3,\alpha_2,\alpha_1,\alpha_0 \right).
\end{align*}
\end{theorem}
We remark that the B{\"a}cklund transformations $u_1,\ldots,u_5$ satisfy the relation \eqref{unvre}. However, the transformations $u_0,u_6$ do not satisfy so.

Theorem \ref{3.1} can be cheched by a direct calculation.

\begin{theorem}\label{3.2}
For the system \eqref{1} of type $D_6^{(1)}$, we make the change of parameters and variables
\begin{gather}
\begin{gathered}\label{27}
\frac{(\alpha_0-\alpha_1)}{2}=A_0, \quad \alpha_1=A_1, \quad \alpha_2=A_2, \quad \alpha_4-\beta_4=A_3, \\ \beta_2=A_4, \quad \beta_4=A_5, \quad \frac{(\beta_3-\beta_4)}{2}=A_6,
\end{gathered}\\
\label{28}
X=\frac{1}{x}, \quad Y=-(xy+\alpha_2)x, \quad Z=\frac{1}{z}, \quad W=-(zw+\beta_2)z
\end{gather}
from $\alpha_0,\alpha_1,\alpha_2,\alpha_4,\beta_2,\beta_3,\beta_4,x,y,z,w$ to $A_0,A_1,\dots ,A_6,X,Y,Z,W$. Then the system \eqref{1} can also be written in the new variables $X,Y,Z,W$ and parameters $A_0,\ldots,A_6$ as a Hamiltonian system. This new system tends to the system \eqref{24} with the Hamiltonian \eqref{25}.
\end{theorem}
{\it Proof.} \qquad Notice that
$$
2(A_0+A_1+A_2+A_3+A_4+A_5+A_6)=\alpha_0+\alpha_1+2\alpha_2+2\alpha_3+2\alpha_4+\alpha_5+\alpha_6=1
$$
and the change of variables from $(x,y,z,w)$ to $(X,Y,Z,W)$ is symplectic. Choose $u_i \ (i=0,1,\dots ,6)$ and $\phi$ as
$$
u_0:=\pi_4, \ u_1:=s_1, \ u_2:=s_2, \ u_3:=s_3, \ u_4:=s_4, \ u_5:=s_6, \ u_6:=\pi_3, \ \phi:=\pi_2.
$$
Then the transformations $u_i$ are reflections of the parameters $A_0,A_1,\dots ,A_6$. The transformation group $\tilde{W}(D_7^{(2)})=<u_0,u_1,\dots ,u_6,\phi>$ coincides with the transformations given in Theorem \ref{3.1}. \qed

Finally, let us summarize Sections 3 and 4 in the below figure.

\begin{figure}[h]
\unitlength 0.1in
\begin{picture}(59.47,60.73)(11.66,-62.55)
%
\special{pn 20}%
\special{ar 2684 376 205 184  1.5707963 6.2831853}%
\special{ar 2684 376 205 184  0.0000000 1.5317917}%
%
\special{pn 20}%
\special{ar 2684 1253 205 185  1.5707963 6.2831853}%
\special{ar 2684 1253 205 185  0.0000000 1.5317917}%
%
\special{pn 20}%
\special{ar 3360 830 205 187  1.5707963 6.2831853}%
\special{ar 3360 830 205 187  0.0000000 1.5317917}%
%
\special{pn 20}%
\special{ar 4901 830 205 187  1.5707963 6.2831853}%
\special{ar 4901 830 205 187  0.0000000 1.5317917}%
%
\special{pn 20}%
\special{ar 5662 368 206 186  1.5707963 6.2831853}%
\special{ar 5662 368 206 186  0.0000000 1.5271348}%
%
\special{pn 20}%
\special{ar 5662 1253 206 185  1.5707963 6.2831853}%
\special{ar 5662 1253 206 185  0.0000000 1.5271348}%
%
\special{pn 20}%
\special{pa 2864 452}%
\special{pa 3198 701}%
\special{fp}%
%
\special{pn 20}%
\special{pa 2872 1186}%
\special{pa 3214 961}%
\special{fp}%
%
\special{pn 20}%
\special{ar 4147 838 206 186  1.5707963 6.2831853}%
\special{ar 4147 838 206 186  0.0000000 1.5271348}%
%
\special{pn 20}%
\special{pa 3565 838}%
\special{pa 3925 838}%
\special{fp}%
%
\special{pn 20}%
\special{pa 4353 838}%
\special{pa 4678 838}%
\special{fp}%
%
\special{pn 20}%
\special{pa 5080 715}%
\special{pa 5483 484}%
\special{fp}%
%
\special{pn 20}%
\special{pa 5055 947}%
\special{pa 5465 1193}%
\special{fp}%
%
\put(21.1000,-4.8500){\makebox(0,0)[lb]{}}%
\put(25.2200,-13.2500){\makebox(0,0)[lb]{$x-t$}}%
\put(32.5800,-8.9700){\makebox(0,0)[lb]{$y$}}%
\put(39.7600,-9.0700){\makebox(0,0)[lb]{$x-z$}}%
\put(47.9800,-9.3500){\makebox(0,0)[lb]{$w$}}%
\put(54.9100,-13.3600){\makebox(0,0)[lb]{$z-1$}}%
\put(55.5200,-4.7900){\makebox(0,0)[lb]{$z$}}%
\put(30.8400,-15.1200){\makebox(0,0)[lb]{Dynkin diagram of type $D_6^{(1)}$}}%
\put(24.8700,-4.3900){\makebox(0,0)[lb]{$x-\infty$}}%
%
\special{pn 20}%
\special{pa 3110 1657}%
\special{pa 2590 2330}%
\special{fp}%
\special{sh 1}%
\special{pa 2590 2330}%
\special{pa 2647 2289}%
\special{pa 2623 2288}%
\special{pa 2615 2265}%
\special{pa 2590 2330}%
\special{fp}%
%
\special{pn 20}%
\special{pa 5150 1669}%
\special{pa 5660 2330}%
\special{fp}%
\special{sh 1}%
\special{pa 5660 2330}%
\special{pa 5635 2265}%
\special{pa 5627 2288}%
\special{pa 5603 2289}%
\special{pa 5660 2330}%
\special{fp}%
%
\special{pn 20}%
\special{ar 3493 2547 124 246  1.5707963 6.2831853}%
\special{ar 3493 2547 124 246  0.0000000 1.5304956}%
%
\special{pn 20}%
\special{ar 3493 3720 124 244  1.5707963 6.2831853}%
\special{ar 3493 3720 124 244  0.0000000 1.5304956}%
%
\special{pn 20}%
\special{pa 2228 3169}%
\special{pa 2445 3169}%
\special{fp}%
%
\special{pn 20}%
\special{pa 2703 3169}%
\special{pa 2900 3169}%
\special{fp}%
%
\special{pn 20}%
\special{pa 3142 3007}%
\special{pa 3385 2699}%
\special{fp}%
%
\special{pn 20}%
\special{pa 3127 3316}%
\special{pa 3375 3642}%
\special{fp}%
\put(20.1000,-32.4800){\makebox(0,0)[lb]{$y$}}%
\put(23.6000,-29.2600){\makebox(0,0)[lb]{$xz-1$}}%
\put(29.7200,-33.0100){\makebox(0,0)[lb]{$w$}}%
\put(32.9000,-34.8500){\makebox(0,0)[lb]{$z-1$}}%
\put(13.4700,-41.4400){\makebox(0,0)[lb]{Dynkin diagram of type $B_6^{(1)}$}}%
%
\special{pn 20}%
\special{ar 2101 3157 124 250  1.5707963 6.2831853}%
\special{ar 2101 3157 124 250  0.0000000 1.5304956}%
%
\special{pn 20}%
\special{pa 1808 3171}%
\special{pa 1959 3171}%
\special{fp}%
%
\special{pn 20}%
\special{pa 1552 3236}%
\special{pa 1394 3236}%
\special{fp}%
\special{sh 1}%
\special{pa 1394 3236}%
\special{pa 1461 3256}%
\special{pa 1447 3236}%
\special{pa 1461 3216}%
\special{pa 1394 3236}%
\special{fp}%
%
\special{pn 20}%
\special{ar 2583 3171 124 245  1.5707963 6.2831853}%
\special{ar 2583 3171 124 245  0.0000000 1.5304956}%
%
\special{pn 20}%
\special{ar 1683 3144 124 247  1.5707963 6.2831853}%
\special{ar 1683 3144 124 247  0.0000000 1.5224469}%
\put(34.0000,-26.4400){\makebox(0,0)[lb]{$z$}}%
\put(16.0000,-32.4800){\makebox(0,0)[lb]{$x$}}%
%
\special{pn 20}%
\special{ar 5860 3201 129 211  1.5707963 6.2831853}%
\special{ar 5860 3201 129 211  0.0000000 1.5320560}%
%
\special{pn 20}%
\special{pa 4567 2768}%
\special{pa 4779 3053}%
\special{fp}%
%
\special{pn 20}%
\special{pa 4572 3603}%
\special{pa 4789 3349}%
\special{fp}%
%
\special{pn 20}%
\special{ar 5381 3207 131 211  1.5707963 6.2831853}%
\special{ar 5381 3207 131 211  0.0000000 1.5250268}%
%
\special{pn 20}%
\special{pa 5011 3207}%
\special{pa 5240 3207}%
\special{fp}%
%
\special{pn 20}%
\special{pa 5512 3207}%
\special{pa 5718 3207}%
\special{fp}%
\put(43.5000,-37.6200){\makebox(0,0)[lb]{$x-t$}}%
\put(48.1700,-32.7500){\makebox(0,0)[lb]{$y$}}%
\put(51.6000,-29.9100){\makebox(0,0)[lb]{$xz-1$}}%
\put(57.9400,-33.1900){\makebox(0,0)[lb]{$w$}}%
\put(47.0600,-39.7700){\makebox(0,0)[lb]{Dynkin diagram of type $B_6^{(1)}$}}%
\put(43.2800,-27.5300){\makebox(0,0)[lb]{$x-\infty$}}%
%
\special{pn 20}%
\special{ar 4450 2633 130 213  1.5707963 6.2831853}%
\special{ar 4450 2633 130 213  0.0000000 1.5323537}%
%
\special{pn 20}%
\special{ar 4450 3679 130 213  1.5707963 6.2831853}%
\special{ar 4450 3679 130 213  0.0000000 1.5323537}%
%
\special{pn 20}%
\special{ar 4884 3191 131 211  1.5707963 6.2831853}%
\special{ar 4884 3191 131 211  0.0000000 1.5250268}%
%
\special{pn 20}%
\special{ar 5862 3218 130 213  1.5707963 6.2831853}%
\special{ar 5862 3218 130 213  0.0000000 1.5246752}%
%
\special{pn 20}%
\special{ar 6942 3203 171 209  1.5707963 6.2831853}%
\special{ar 6942 3203 171 209  0.0000000 1.5298835}%
%
\special{pn 20}%
\special{pa 5994 3229}%
\special{pa 6150 3229}%
\special{fp}%
%
\special{pn 20}%
\special{pa 6523 3161}%
\special{pa 6769 3161}%
\special{fp}%
\special{sh 1}%
\special{pa 6769 3161}%
\special{pa 6702 3141}%
\special{pa 6716 3161}%
\special{pa 6702 3181}%
\special{pa 6769 3161}%
\special{fp}%
%
\special{pn 20}%
\special{pa 6515 3294}%
\special{pa 6777 3294}%
\special{fp}%
\special{sh 1}%
\special{pa 6777 3294}%
\special{pa 6710 3274}%
\special{pa 6724 3294}%
\special{pa 6710 3314}%
\special{pa 6777 3294}%
\special{fp}%
\put(61.5000,-29.8800){\makebox(0,0)[lb]{$z-\infty$}}%
%
\special{pn 20}%
\special{pa 2580 4417}%
\special{pa 3100 5273}%
\special{fp}%
\special{sh 1}%
\special{pa 3100 5273}%
\special{pa 3082 5206}%
\special{pa 3072 5227}%
\special{pa 3048 5226}%
\special{pa 3100 5273}%
\special{fp}%
%
\special{pn 20}%
\special{pa 5660 4430}%
\special{pa 5150 5325}%
\special{fp}%
\special{sh 1}%
\special{pa 5150 5325}%
\special{pa 5200 5277}%
\special{pa 5176 5279}%
\special{pa 5166 5257}%
\special{pa 5150 5325}%
\special{fp}%
%
\special{pn 20}%
\special{ar 3343 5742 205 187  1.5707963 6.2831853}%
\special{ar 3343 5742 205 187  0.0000000 1.5317917}%
%
\special{pn 20}%
\special{ar 4884 5742 205 187  1.5707963 6.2831853}%
\special{ar 4884 5742 205 187  0.0000000 1.5317917}%
%
\special{pn 20}%
\special{ar 4130 5749 206 188  1.5707963 6.2831853}%
\special{ar 4130 5749 206 188  0.0000000 1.5271348}%
%
\special{pn 20}%
\special{pa 3548 5749}%
\special{pa 3908 5749}%
\special{fp}%
%
\special{pn 20}%
\special{pa 4336 5749}%
\special{pa 4661 5749}%
\special{fp}%
\put(32.4100,-58.1000){\makebox(0,0)[lb]{$y$}}%
\put(39.5900,-58.2100){\makebox(0,0)[lb]{$x-z$}}%
\put(47.8100,-58.4800){\makebox(0,0)[lb]{$w$}}%
\put(30.6700,-64.2500){\makebox(0,0)[lb]{Dynkin diagram of type $D_7^{(2)}$}}%
%
\special{pn 20}%
\special{ar 2610 5740 205 186  1.5707963 6.2831853}%
\special{ar 2610 5740 205 186  0.0000000 1.5317917}%
%
\special{pn 20}%
\special{ar 5660 5753 205 185  1.5707963 6.2831853}%
\special{ar 5660 5753 205 185  0.0000000 1.5317917}%
%
\special{pn 20}%
\special{ar 1800 5753 205 185  1.5707963 6.2831853}%
\special{ar 1800 5753 205 185  0.0000000 1.5317917}%
%
\special{pn 20}%
\special{ar 6450 5740 205 186  1.5707963 6.2831853}%
\special{ar 6450 5740 205 186  0.0000000 1.5317917}%
%
\special{pn 20}%
\special{pa 2400 5687}%
\special{pa 2000 5687}%
\special{fp}%
\special{sh 1}%
\special{pa 2000 5687}%
\special{pa 2067 5707}%
\special{pa 2053 5687}%
\special{pa 2067 5667}%
\special{pa 2000 5687}%
\special{fp}%
%
\special{pn 20}%
\special{pa 5850 5869}%
\special{pa 6260 5869}%
\special{fp}%
\special{sh 1}%
\special{pa 6260 5869}%
\special{pa 6193 5849}%
\special{pa 6207 5869}%
\special{pa 6193 5889}%
\special{pa 6260 5869}%
\special{fp}%
%
\special{pn 20}%
\special{pa 2830 5765}%
\special{pa 3110 5765}%
\special{fp}%
%
\special{pn 20}%
\special{pa 5090 5765}%
\special{pa 5420 5765}%
\special{fp}%
\put(24.7000,-58.4400){\makebox(0,0)[lb]{$x$}}%
\put(53.9000,-55.6900){\makebox(0,0)[lb]{$z-\infty$}}%
%
\special{pn 20}%
\special{ar 6341 3219 171 210  1.5707963 6.2831853}%
\special{ar 6341 3219 171 210  0.0000000 1.5298835}%
%
\special{pn 20}%
\special{ar 3030 3173 124 245  1.5707963 6.2831853}%
\special{ar 3030 3173 124 245  0.0000000 1.5304956}%
%
\special{pn 20}%
\special{ar 1290 3121 124 244  1.5707963 6.2831853}%
\special{ar 1290 3121 124 244  0.0000000 1.5304956}%
%
\special{pn 20}%
\special{pa 2400 5818}%
\special{pa 2010 5818}%
\special{fp}%
\special{sh 1}%
\special{pa 2010 5818}%
\special{pa 2077 5838}%
\special{pa 2063 5818}%
\special{pa 2077 5798}%
\special{pa 2010 5818}%
\special{fp}%
%
\special{pn 20}%
\special{pa 5880 5722}%
\special{pa 6250 5722}%
\special{fp}%
\special{sh 1}%
\special{pa 6250 5722}%
\special{pa 6183 5702}%
\special{pa 6197 5722}%
\special{pa 6183 5742}%
\special{pa 6250 5722}%
\special{fp}%
%
\special{pn 20}%
\special{pa 3970 1657}%
\special{pa 3970 5299}%
\special{fp}%
\special{sh 1}%
\special{pa 3970 5299}%
\special{pa 3990 5232}%
\special{pa 3970 5246}%
\special{pa 3950 5232}%
\special{pa 3970 5299}%
\special{fp}%
%
\special{pn 20}%
\special{pa 1550 3053}%
\special{pa 1400 3053}%
\special{fp}%
\special{sh 1}%
\special{pa 1400 3053}%
\special{pa 1467 3073}%
\special{pa 1453 3053}%
\special{pa 1467 3033}%
\special{pa 1400 3053}%
\special{fp}%
%
\special{pn 8}%
\special{pa 1820 838}%
\special{pa 3120 838}%
\special{dt 0.045}%
\special{pa 3120 838}%
\special{pa 3119 838}%
\special{dt 0.045}%
%
\special{pn 8}%
\special{pa 5110 861}%
\special{pa 6400 861}%
\special{dt 0.045}%
\special{pa 6400 861}%
\special{pa 6399 861}%
\special{dt 0.045}%
%
\special{pn 20}%
\special{pa 2210 676}%
\special{pa 2189 651}%
\special{pa 2166 629}%
\special{pa 2139 612}%
\special{pa 2106 602}%
\special{pa 2072 601}%
\special{pa 2042 611}%
\special{pa 2020 631}%
\special{pa 2005 659}%
\special{pa 1995 692}%
\special{pa 1989 726}%
\special{pa 1983 759}%
\special{pa 1980 791}%
\special{pa 1978 823}%
\special{pa 1979 854}%
\special{pa 1983 885}%
\special{pa 1991 917}%
\special{pa 2002 950}%
\special{pa 2017 983}%
\special{pa 2037 1013}%
\special{pa 2060 1035}%
\special{pa 2089 1044}%
\special{pa 2121 1037}%
\special{pa 2152 1018}%
\special{pa 2174 991}%
\special{pa 2180 959}%
\special{pa 2180 956}%
\special{sp}%
%
\special{pn 20}%
\special{pa 6120 687}%
\special{pa 6139 655}%
\special{pa 6159 626}%
\special{pa 6180 603}%
\special{pa 6202 590}%
\special{pa 6226 589}%
\special{pa 6252 598}%
\special{pa 6277 617}%
\special{pa 6301 643}%
\special{pa 6323 676}%
\special{pa 6342 713}%
\special{pa 6357 753}%
\special{pa 6367 794}%
\special{pa 6372 836}%
\special{pa 6372 877}%
\special{pa 6367 915}%
\special{pa 6358 949}%
\special{pa 6344 978}%
\special{pa 6327 1001}%
\special{pa 6305 1017}%
\special{pa 6280 1023}%
\special{pa 6252 1021}%
\special{pa 6222 1012}%
\special{pa 6190 999}%
\special{pa 6157 982}%
\special{pa 6150 979}%
\special{sp}%
%
\special{pn 20}%
\special{pa 2180 979}%
\special{pa 2230 827}%
\special{fp}%
\special{sh 1}%
\special{pa 2230 827}%
\special{pa 2190 884}%
\special{pa 2213 878}%
\special{pa 2228 897}%
\special{pa 2230 827}%
\special{fp}%
%
\special{pn 20}%
\special{pa 6160 990}%
\special{pa 6050 945}%
\special{fp}%
\special{sh 1}%
\special{pa 6050 945}%
\special{pa 6104 989}%
\special{pa 6099 965}%
\special{pa 6119 952}%
\special{pa 6050 945}%
\special{fp}%
\put(12.2000,-32.2500){\makebox(0,0)[lb]{$\pi_4$}}%
\put(68.3000,-33.1200){\makebox(0,0)[lb]{$\pi_3$}}%
\put(16.6000,-58.5100){\makebox(0,0)[lb]{$\pi_4$}}%
\put(63.1000,-58.5100){\makebox(0,0)[lb]{$\pi_3$}}%
\put(18.1000,-6.0100){\makebox(0,0)[lb]{$\pi_4$}}%
\put(62.8000,-6.2200){\makebox(0,0)[lb]{$\pi_3$}}%
\end{picture}%
\label{fig:CPVI13}
\caption{Dynkin diagrams of types $B_6^{(1)}$,$D_6^{(1)}$ and $D_7^{(2)}$}
\end{figure}
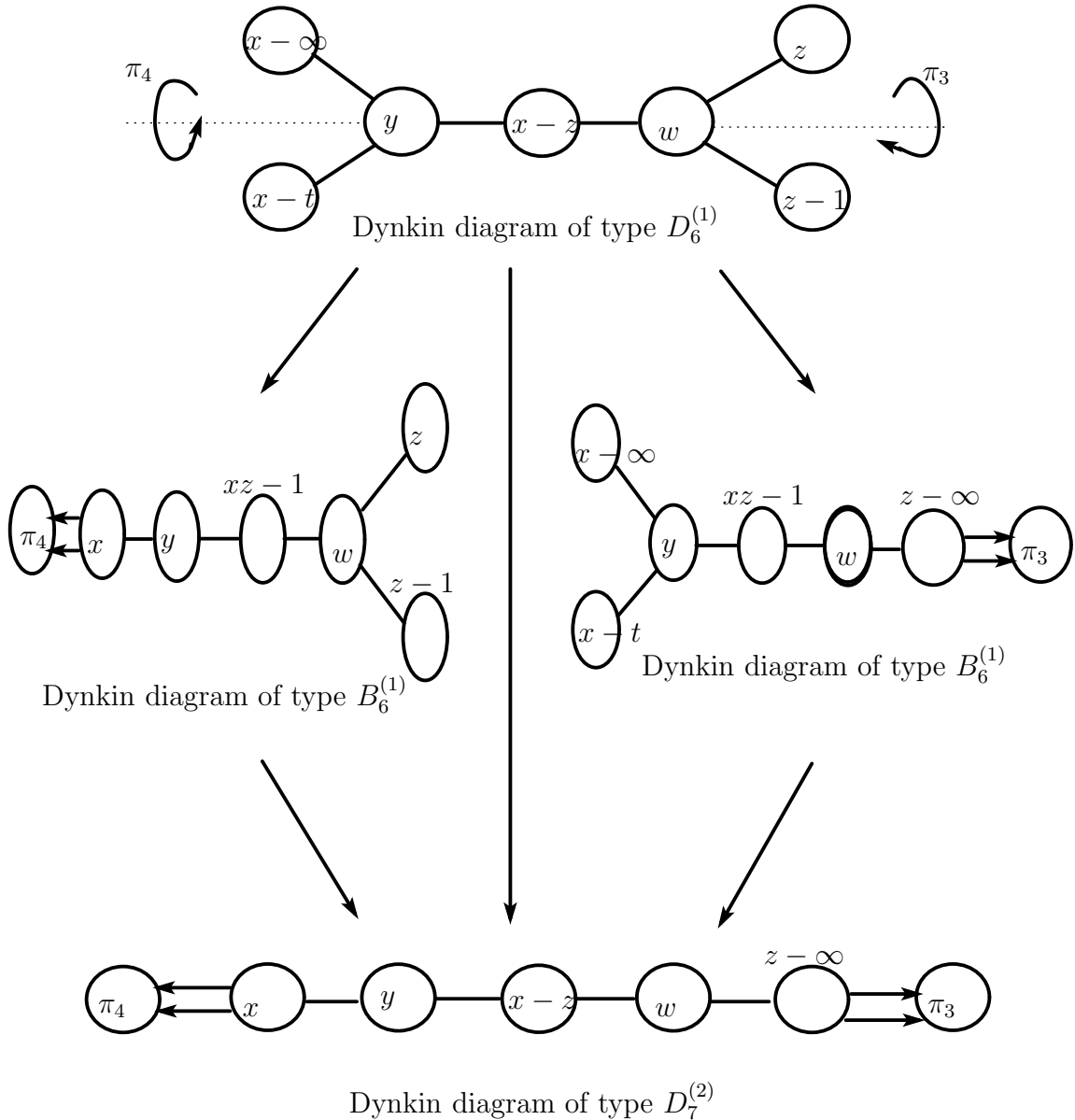

\section{Autonomous version of type $D_6^{(1)}$}

In this section, we find an autonomous version of the system \eqref{1} explicitly given by
\begin{equation}\label{t1}
\frac{dx}{dt}=\frac{\partial H}{\partial y}, \ \ \frac{dy}{dt}=-\frac{\partial H}{\partial x}, \ \ \frac{dz}{dt}=\frac{\partial H}{\partial w}, \ \ \frac{dw}{dt}=-\frac{\partial H}{\partial z}
\end{equation}
with the polynomial Hamiltonian
\begin{align}\label{t2}
\begin{split}
H &=H_{auto}(x,y;\alpha_0,\alpha_1,\alpha_2,\alpha_3+2\alpha_4+\alpha_5,\alpha_3+\alpha_6)\\
&+H_{auto}(z,w;\alpha_0+\alpha_3,\alpha_1+2\alpha_2+\alpha_3,\alpha_4,\alpha_5,\alpha_6)\\
&+2(x-\eta)yz\{(z-1)w+\alpha_4 \}.
  \end{split}
\end{align}
Here $x,y,z$ and $w$ denote unknown complex variables, and $\eta$ and $\alpha_0,\alpha_1,\dots ,\alpha_6$ are complex parameters satisfying the relation:
\begin{equation}\label{t3}
\begin{aligned}
\alpha_0+\alpha_1+2\alpha_2+2\alpha_3+2\alpha_4+\alpha_5+\alpha_6=0.
\end{aligned}
\end{equation}
Here the symbol $H_{auto}(q,p;\beta_0,\beta_1,\beta_2,\beta_3,\beta_4)$ denotes the Hamiltonian given by
\begin{equation}\label{t4}
\begin{aligned}
&H_{auto}(q,p;\beta_0,\beta_1,\beta_2,\beta_3,\beta_4)=\\
&q(q-1)(q-\eta)p^2-\{\beta_0q(q-1)+\beta_4(q-1)(q-\eta)+\beta_3q(q-\eta)\}p+\beta_2(\beta_1+\beta_2)q\\
&\hspace{9cm} (\beta_0+\beta_1+2\beta_2+\beta_3+\beta_4=0).
\end{aligned}
\end{equation}

\begin{proposition}\label{t1.3}
The system \eqref{t1} has the Hamiltonian \eqref{t2} as its first integral.
\end{proposition}

\begin{theorem}\label{t1.1}
The system \eqref{t1} admits the affine Weyl group symmetry of type $D_6^{(1)}$ as the group of its B{\"a}cklund transformations, whose generators $s_0,s_1,\dots ,s_6$ are explicitly given as follows{\rm : \rm}with the notation$:(*):=(x,y,z,w;\alpha_0,\alpha_1,\ldots,\alpha_6),$
\begin{align*}
s_0: (*) &\rightarrow \left(x,y-\frac{\alpha_0}{x-\eta},z,w;-\alpha_0,\alpha_1,\alpha_2+\alpha_0,\alpha_3,\alpha_4,\alpha_5,\alpha_6 \right), \\
s_1: (*) &\rightarrow (x,y,z,w;\alpha_0,-\alpha_1,\alpha_2+\alpha_1,\alpha_3,\alpha_4,\alpha_5,\alpha_6), \\
s_2: (*) &\rightarrow \left(x+\frac{\alpha_2}{y},y,z,w;\alpha_0+\alpha_2,\alpha_1+\alpha_2,-\alpha_2,\alpha_3+\alpha_2,\alpha_4,\alpha_5,\alpha_6 \right), \\
s_3: (*) &\rightarrow \left(x,y-\frac{\alpha_3}{x-z},z,w+\frac{\alpha_3}{x-z};\alpha_0,\alpha_1,\alpha_2+\alpha_3,-\alpha_3,\alpha_4+\alpha_3,\alpha_5,\alpha_6 \right), \\
s_4: (*) &\rightarrow \left(x,y,z+\frac{\alpha_4}{w},w;\alpha_0,\alpha_1,\alpha_2,\alpha_3+\alpha_4,-\alpha_4,\alpha_5+\alpha_4,\alpha_6+\alpha_4 \right), \\
s_5: (*) &\rightarrow \left(x,y,z,w-\frac{\alpha_5}{z-1};\alpha_0,\alpha_1,\alpha_2,\alpha_3,\alpha_4+\alpha_5,-\alpha_5,\alpha_6 \right), \\
s_6: (*) &\rightarrow \left(x,y,z,w-\frac{\alpha_6}{z};\alpha_0,\alpha_1,\alpha_2,\alpha_3,\alpha_4+\alpha_6,\alpha_5,-\alpha_6 \right). 
\end{align*}
\end{theorem}
\noindent
We note that these transformations $s_i$ are birational and symplectic.

\begin{theorem}\label{t1.2}
Let us consider a polynomial Hamiltonian system with Hamiltonian $H \in {\Bbb C}[x,y,z,w]$. We assume that

$(A1)$ $deg(H)=5$ with respect to $x,y,z,w$.

$(A2)$ This system becomes again a polynomial Hamiltonian system in each coordinate system $(x_i,y_i,z_i,w_i) \ (i=0,1,\dots ,6)${\rm : \rm}
\begin{align*}
&x_0=-((x-\eta)y-\alpha_0)y, \quad y_0=1/y, \quad z_0=z, \quad w_0=w, \\
&x_1=1/x, \quad y_1=-x(xy+\alpha_1+\alpha_2), \quad z_1=z, \quad w_1=w, \\
&x_2=1/x, \quad y_2=-x(xy+\alpha_2), \quad z_2=z, \quad w_2=w, \\
&x_3=-((x-z)y-\alpha_3)y, \quad y_3=1/y, \quad z_3=z, \quad w_3=y+w, \\
&x_4=x, \quad y_4=y, \quad z_4=1/z, \quad w_4=-z(zw+\alpha_4), \\
&x_5=x, \quad y_5=y, \quad z_5=-((z-1)w-\alpha_5)w, \quad w_5=1/w, \\
&x_6=x, \quad y_6=y, \quad z_6=-w(zw-\alpha_6), \quad w_6=1/w.
\end{align*}
Then such a system coincides with the system \eqref{t1}.
\end{theorem}

\section{Review of the systems of types $A_4^{(1)}$ and $A_5^{(1)}$}

Let us recall the system of type $A_5^{(1)}$ given by
\begin{equation}\label{29}
  \left\{
  \begin{aligned}
   \frac{dx}{dt} &=\frac{\partial H_{A_5^{(1)}}}{\partial y}=\frac{2x^2y+2xzw}{t}-\frac{x^2}{t}-2xy-2zw+\left(1+\frac{\alpha_1+\alpha_3+\alpha_5}{t} \right)x+\alpha_2+\alpha_4,\\
   \frac{dy}{dt} &=-\frac{\partial H_{A_5^{(1)}}}{\partial x}=\frac{-2xy^2-2yzw}{t}+y^2+\frac{2xy}{t}-\left(1+\frac{\alpha_1+\alpha_3+\alpha_5}{t} \right)y+\frac{\alpha_1}{t},\\
   \frac{dz}{dt} &=\frac{\partial H_{A_5^{(1)}}}{\partial w}=\frac{2z^2w+2xyz}{t}-\frac{z^2}{t}-2zw-2yz+\left(1+\frac{\alpha_1+\alpha_3+\alpha_5}{t} \right)z+\alpha_4,\\
   \frac{dw}{dt} &=-\frac{\partial H_{A_5^{(1)}}}{\partial z}=\frac{-2zw^2-2xyw}{t}+w^2+\frac{2zw}{t}+2yw-\left(1+\frac{\alpha_1+\alpha_3+\alpha_5}{t} \right)w+\frac{\alpha_3}{t}
   \end{aligned}
  \right. 
\end{equation}
with the polynomial Hamiltonian:
\begin{align}\label{30}
\begin{split}
H_{A_5^{(1)}}=&H_V(x,y,t;\alpha_1+\alpha_3+\alpha_5,\alpha_2+\alpha_4,\alpha_1)\\
&+H_V(z,w,t;\alpha_1+\alpha_3+\alpha_5,\alpha_4,\alpha_3)-2yzw+\frac{2xyzw}{t}.
\end{split}
\end{align}
Here, $x,y,z$ and $w$ denote unknown complex variables, and $\alpha_0,\alpha_1,\dots ,\alpha_5$ are complex parameters with $\alpha_0+\alpha_1+\alpha_2+\alpha_3+\alpha_4+\alpha_5=1$, and the symbol $H_V(q,p,t;\gamma_1,\gamma_2,\gamma_3)$ denotes the Hamiltonian of the second-order Painlev\'e V equation given by
\begin{equation}\label{31}
\begin{aligned}
H_{V}(q,p,t;\gamma_1,\gamma_2,\gamma_3)=\frac{q^2p^2-q^2p}{t}-qp^2+\left(1+\frac{\gamma_3}{t} \right)qp+\gamma_2p-\frac{\gamma_1q}{t}.
\end{aligned}
\end{equation}

\begin{figure}[h]
\unitlength 0.1in
\begin{picture}( 33.5100, 19.5000)( 31.3900,-23.9000)
%
\special{pn 20}%
\special{ar 4868 810 330 158  1.6079080 6.2831853}%
\special{ar 4868 810 330 158  0.0000000 1.6066661}%
%
\special{pn 20}%
\special{ar 3470 1230 332 158  1.6076347 6.2831853}%
\special{ar 3470 1230 332 158  0.0000000 1.6063264}%
%
\special{pn 20}%
\special{ar 3470 1828 332 158  1.6076347 6.2831853}%
\special{ar 3470 1828 332 158  0.0000000 1.6063264}%
%
\special{pn 20}%
\special{ar 4854 2234 330 158  1.6077253 6.2831853}%
\special{ar 4854 2234 330 158  0.0000000 1.6064954}%
%
\special{pn 20}%
\special{ar 6160 1216 330 158  1.6103796 6.2831853}%
\special{ar 6160 1216 330 158  0.0000000 1.6090552}%
%
\special{pn 20}%
\special{ar 6160 1814 330 158  1.6103796 6.2831853}%
\special{ar 6160 1814 330 158  0.0000000 1.6088732}%
%
\special{pn 20}%
\special{pa 4582 896}%
\special{pa 3696 1108}%
\special{fp}%
%
\special{pn 20}%
\special{pa 3682 1948}%
\special{pa 4522 2170}%
\special{fp}%
%
\special{pn 20}%
\special{pa 5124 902}%
\special{pa 5904 1116}%
\special{fp}%
%
\special{pn 20}%
\special{pa 6144 1380}%
\special{pa 6144 1642}%
\special{fp}%
%
\special{pn 20}%
\special{pa 3440 1394}%
\special{pa 3440 1664}%
\special{fp}%
%
\special{pn 20}%
\special{pa 6010 1956}%
\special{pa 5184 2176}%
\special{fp}%
\put(32.0000,-12.9300){\makebox(0,0)[lb]{$x-t$}}%
\put(33.2000,-18.9200){\makebox(0,0)[lb]{$y$}}%
\put(46.1000,-23.0000){\makebox(0,0)[lb]{$x-z$}}%
\put(60.3900,-18.7700){\makebox(0,0)[lb]{$w$}}%
\put(45.4000,-6.1000){\makebox(0,0)[lb]{$y+w-1$}}%
\put(60.0900,-12.7200){\makebox(0,0)[lb]{$z$}}%
\end{picture}%
\label{fig:CPVI18}
\caption{Dynkin diagram of type $A_5^{(1)}$}
\end{figure}
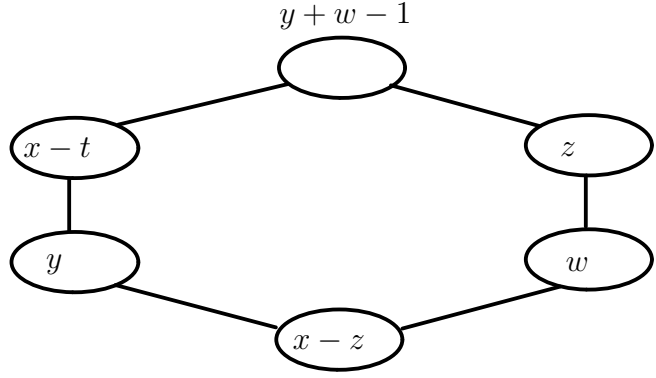

It is known that the system \eqref{29} admits the affine Weyl group symmetry of type $A_5^{(1)}$ as the group of its B{\"a}cklund transformations, whose generators $s_0,s_1,\dots ,s_5$ are explicitly given as follows{\rm : \rm}with the notation$:(*):=(x,y,z,w,t;\alpha_0,\alpha_1,\ldots,\alpha_5),$
\begin{align*}
s_0: (*) &\rightarrow \left(x,y-\frac{\alpha_0}{x-t},z,w,t;-\alpha_0,\alpha_1+\alpha_0,\alpha_2,\alpha_3,\alpha_4,\alpha_5+\alpha_0 \right), \\
s_1: (*) &\rightarrow \left(x+\frac{\alpha_1}{y},y,z,w,t;\alpha_0+\alpha_1,-\alpha_1,\alpha_2+\alpha_1,\alpha_3,\alpha_4,\alpha_5 \right), \\
s_2: (*) &\rightarrow \left(x,y-\frac{\alpha_2}{x-z},z,w+\frac{\alpha_2}{x-z},t;\alpha_0,\alpha_1+\alpha_2,-\alpha_2,\alpha_3+\alpha_2,\alpha_4,\alpha_5 \right), \\
s_3: (*) &\rightarrow \left(x,y,z+\frac{\alpha_3}{w},w,t;\alpha_0,\alpha_1,\alpha_2+\alpha_3,-\alpha_3,\alpha_4+\alpha_3,\alpha_5 \right), \\
s_4: (*) &\rightarrow \left(x,y,z,w-\frac{\alpha_4}{z},t;\alpha_0,\alpha_1,\alpha_2,\alpha_3+\alpha_4,-\alpha_4,\alpha_5+\alpha_4 \right), \\
s_5: (*) &\rightarrow \left(x+\frac{\alpha_5}{y+w-1},y,z+\frac{\alpha_5}{y+w-1},w,t;\alpha_0+\alpha_5,\alpha_1,\alpha_2,\alpha_3,\alpha_4+\alpha_5,-\alpha_5 \right).
\end{align*}

\begin{theorem}\label{4.1}
Let us consider a polynomial Hamiltonian system with Hamiltonian $H \in C(t)[x,y,z,w]$. We assume that

$(A1)$ $deg(H)=4$ with respect to $x,y,z,w$.

$(A2)$ This system becomes again a polynomial Hamiltonian system in each coordinate system $r_i(i=0,1,\dots ,5)${\rm : \rm}
\begin{align*}
&r_0:x_0=-((x-t)y-\alpha_0)y, \quad y_0=1/y, \quad z_0=z, \quad w_0=w,\\
&r_1:x_1=1/x, \quad y_1=-(xy+\alpha_1)x, \quad z_1=z, \quad w_1=w,\\
&r_2:x_2=-((x-z)y-\alpha_2)y, \quad y_2=1/y, \quad z_2=z, \quad w_2=w+y,\\
&r_3:x_3=x, \quad y_3=y, \quad z_3=1/z, \quad w_3=-(zw+\alpha_3)z,\\
&r_4:x_4=x, \quad y_4=y, \quad z_4=-(zw-\alpha_4)w, \quad w_4=1/w,\\
&r_5:x_5=1/x, \quad y_5=-((y+w-1)y+\alpha_5)x, \quad z_5=z-x, \quad w_5=w.
\end{align*}
Then such a system coincides with the system \eqref{29}.
\end{theorem}

Next, let us recall the system of type $A_4^{(1)}$ given by
\begin{equation}\label{32}
  \left\{
  \begin{aligned}
   \frac{dx}{dt} &=\frac{\partial H_{A_4^{(1)}}}{\partial y}=x^2+2xy+2zw-tx+\alpha_2+\alpha_4,\\
   \frac{dy}{dt} &=-\frac{\partial H_{A_4^{(1)}}}{\partial x}=-y^2-2xy+ty+\alpha_1,\\
   \frac{dz}{dt} &=\frac{\partial H_{A_4^{(1)}}}{\partial w}=z^2+2zw+2yz-tz+\alpha_4,\\
   \frac{dw}{dt} &=-\frac{\partial H_{A_4^{(1)}}}{\partial z}=-w^2-2zw-2yw+tw+\alpha_3
   \end{aligned}
  \right. 
\end{equation}
with the polynomial Hamiltonian:
\begin{equation}\label{33}
\begin{aligned}
H_{A_4^{(1)}}=H_{IV}(x,y,t;\alpha_1,\alpha_2+\alpha_4)+H_{IV}(z,w,t;\alpha_3,\alpha_4)+2yzw.\\
\end{aligned}
\end{equation}
Here, $x,y,z$ and $w$ denote unknown complex variables, and $\alpha_0,\alpha_1,\dots ,\alpha_4$ are complex parameters with $\alpha_0+\alpha_1+\alpha_2+\alpha_3+\alpha_4=1$, and the symbol $H_{IV}(q,p,t;\gamma_1,\gamma_2)$ denotes the Hamiltonian of the second-order Painlev\'e IV equation given by
\begin{equation}\label{34}
\begin{aligned}
H_{IV}(q,p,t;\gamma_1,\gamma_2)=q^2p+qp^2-tqp-\gamma_1q+\gamma_2p.
\end{aligned}
\end{equation}
It is known that the system \eqref{32} admits the affine Weyl group symmetry of type $A_4^{(1)}$ as the group of its B{\"a}cklund transformations, whose generators $s_0,s_1,\dots ,s_4$ are explicitly given as follows{\rm : \rm}with the notation$:(*):=(x,y,z,w,t;\alpha_0,\alpha_1,\ldots,\alpha_4),$
\begin{align*}
\begin{split}
s_0: (*) &\rightarrow (x+\frac{\alpha_0}{x+y+w-t},y-\frac{\alpha_0}{x+y+w-t},z+\frac{\alpha_0}{x+y+w-t},w,t; \\
& \qquad -\alpha_0,\alpha_1+\alpha_0,\alpha_2,\alpha_3,\alpha_4+\alpha_0),
\end{split}\\
s_1: (*) &\rightarrow \left(x+\frac{\alpha_1}{y},y,z,w,t;\alpha_0+\alpha_1,-\alpha_1,\alpha_2+\alpha_1,\alpha_3,\alpha_4 \right), \\
s_2: (*) &\rightarrow \left(x,y-\frac{\alpha_2}{x-z},z,w+\frac{\alpha_2}{x-z},t;\alpha_0,\alpha_1+\alpha_2,-\alpha_2,\alpha_3+\alpha_2,\alpha_4 \right), \\
s_3: (*) &\rightarrow \left(x,y,z+\frac{\alpha_3}{w},w,t;\alpha_0,\alpha_1,\alpha_2+\alpha_3,-\alpha_3,\alpha_4+\alpha_3 \right), \\
s_4: (*) &\rightarrow \left(x,y,z,w-\frac{\alpha_4}{z},t;\alpha_0+\alpha_4,\alpha_1,\alpha_2,\alpha_3+\alpha_4,-\alpha_4 \right).
\end{align*}

\begin{figure}
\unitlength 0.1in
\begin{picture}(23.60,20.50)(32.60,-25.00)
%
\special{pn 20}%
\special{ar 3610 870 220 220  1.6078164 6.2831853}%
\special{ar 3610 870 220 220  0.0000000 1.6064954}%
%
\special{pn 20}%
\special{ar 3610 1710 220 220  1.6078164 6.2831853}%
\special{ar 3610 1710 220 220  0.0000000 1.6064954}%
%
\special{pn 20}%
\special{ar 4530 2280 220 220  1.6078164 6.2831853}%
\special{ar 4530 2280 220 220  0.0000000 1.6064954}%
%
\special{pn 20}%
\special{ar 5400 850 220 220  1.6078164 6.2831853}%
\special{ar 5400 850 220 220  0.0000000 1.6064954}%
%
\special{pn 20}%
\special{ar 5400 1690 220 220  1.6078164 6.2831853}%
\special{ar 5400 1690 220 220  0.0000000 1.6064954}%
%
\special{pn 20}%
\special{pa 3750 1880}%
\special{pa 4310 2190}%
\special{fp}%
%
\special{pn 20}%
\special{pa 5390 1080}%
\special{pa 5390 1450}%
\special{fp}%
%
\special{pn 20}%
\special{pa 3590 1100}%
\special{pa 3590 1480}%
\special{fp}%
%
\special{pn 20}%
\special{pa 5300 1890}%
\special{pa 4750 2200}%
\special{fp}%
\put(32.6000,-6.2000){\makebox(0,0)[lb]{$x+y+w-t$}}%
\put(35.1000,-18.0000){\makebox(0,0)[lb]{$y$}}%
\put(43.4000,-23.6000){\makebox(0,0)[lb]{$x-z$}}%
\put(53.2000,-17.8000){\makebox(0,0)[lb]{$w$}}%
\put(53.0000,-9.3000){\makebox(0,0)[lb]{$z$}}%
%
\special{pn 20}%
\special{pa 3830 860}%
\special{pa 5170 860}%
\special{fp}%
\end{picture}%
\label{fig:CPVI19}
\caption{Dynkin diagram of type $A_4^{(1)}$}
\end{figure}
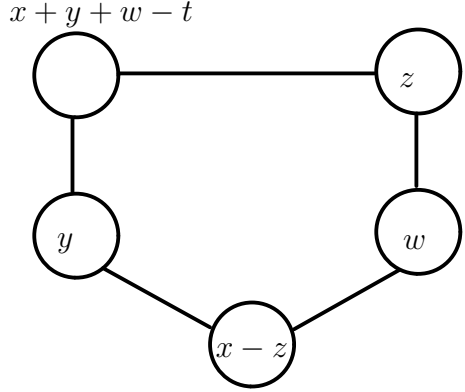

\begin{theorem}\label{4.2}
Let us consider a polynomial Hamiltonian system with Hamiltonian $H \in C(t)[x,y,z,w]$. We assume that

$(A1)$ $deg(H)=3$ with respect to $x,y,z,w$.

$(A2)$ This system becomes again a polynomial Hamiltonian system in each coordinate system $r_i(i=0,1,\dots ,4)${\rm : \rm}
\begin{align*}
&r_0:x_0=-((x+y+w-t)y-\alpha_0)y, \quad y_0=1/y, \quad z_0=z+y, \quad w_0=w,\\
&r_1:x_1=1/x, \quad y_1=-(xy+\alpha_1)x, \quad z_1=z, \quad w_1=w,\\
&r_2:x_2=-((x-z)y-\alpha_2)y, \quad y_2=1/y, \quad z_2=z, \quad w_2=w+y,\\
&r_3:x_3=x, \quad y_3=y, \quad z_3=1/z, \quad w_3=-(zw+\alpha_3)z,\\
&r_4:x_4=x, \quad y_4=y, \quad z_4=-(zw-\alpha_4)w, \quad w_4=1/w.
\end{align*}
Then such a system coincides with the system \eqref{32}.
\end{theorem}
Theorems \ref{4.1} and \ref{4.2} can be cheched by a direct calculation, respectively.

\section{An approach for obtaining the system \eqref{1}}

Much effort has been made to investigate algebraic ordinary differential systems with symmetry under the affine Weyl group of type $D_6^{(1)}$, however these systems have not yet been found. Taking a hint from the representation of the affine Weyl groups of types $A_4^{(1)}$ and $A_5^{(1)}$, we consider Problem 1. We do not yet have explicit descriptions of the symmetry under the affine Weyl group of type $D_6^{(1)}$ with respect to $x,y,z,w$, so we will construct a representation under the affine Weyl group of type $D_6^{(1)}$ by using a part of the symmetry under the affine Weyl groups of types $A_4^{(1)}$ and $A_5^{(1)}$. In the case of the Painlev\'e systems, the fourth, fifth and sixth Painlev\'e systems have affine Weyl group symmetry of type $A_2^{(1)},A_3^{(1)}$ and $D_4^{(1)}$, respectively. Each of them has a common subgroup, which is isomorphic to the classical Weyl group $W(A_2)$. Here, the elements $u_i$ of the subgroup $W(A_2)=<u_1,u_2>$ are explicitly written as follows:
\begin{align*}
\begin{split}
u_1: (x,y,\gamma_1,\gamma_2) &\rightarrow \left(x+\frac{\gamma_1}{y},y,-\gamma_1,\gamma_2+\gamma_1 \right), \\ u_2:(x,y,\gamma_1,\gamma_2) &\rightarrow \left(x,y-\frac{\gamma_2}{x},\gamma_1+\gamma_2,-\gamma_2 \right).
\end{split}
\end{align*}
Here, $\gamma_1$ and $\gamma_2$ are root parameters.

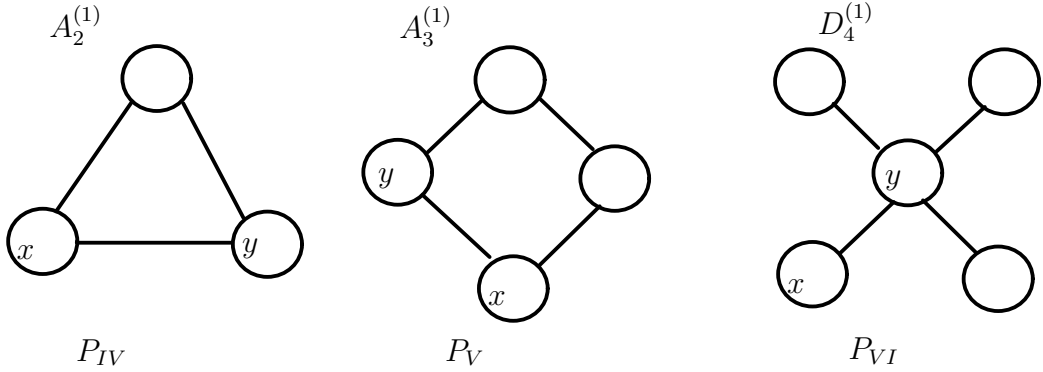
\begin{figure}[h]
\unitlength 0.1in
\begin{picture}(53.96,16.95)(6.44,-19.55)
%
\special{pn 20}%
\special{ar 4838 641 180 170  0.0000000 6.2831853}%
%
\special{pn 20}%
\special{ar 5353 1118 180 170  0.0000000 6.2831853}%
%
\special{pn 20}%
\special{ar 5859 641 181 170  0.0000000 6.2831853}%
%
\special{pn 20}%
\special{ar 4854 1649 180 169  0.0000000 6.2831853}%
%
\special{pn 20}%
\special{ar 5826 1673 181 169  0.0000000 6.2831853}%
%
\special{pn 20}%
\special{ar 3271 631 181 171  0.0000000 6.2831853}%
%
\special{pn 20}%
\special{ar 2683 1115 180 170  0.0000000 6.2831853}%
%
\special{pn 20}%
\special{ar 3288 1724 180 169  0.0000000 6.2831853}%
%
\special{pn 20}%
\special{ar 3819 1148 180 169  0.0000000 6.2831853}%
%
\special{pn 20}%
\special{ar 1422 624 180 170  0.0000000 6.2831853}%
%
\special{pn 20}%
\special{ar 825 1477 181 171  0.0000000 6.2831853}%
%
\special{pn 20}%
\special{ar 2001 1492 181 171  0.0000000 6.2831853}%
%
\special{pn 20}%
\special{pa 4964 756}%
\special{pa 5198 988}%
\special{fp}%
%
\special{pn 20}%
\special{pa 5756 773}%
\special{pa 5504 1005}%
\special{fp}%
%
\special{pn 20}%
\special{pa 5000 1537}%
\special{pa 5279 1271}%
\special{fp}%
%
\special{pn 20}%
\special{pa 5711 1537}%
\special{pa 5432 1262}%
\special{fp}%
%
\special{pn 20}%
\special{pa 3128 735}%
\special{pa 2840 1010}%
\special{fp}%
%
\special{pn 20}%
\special{pa 2813 1234}%
\special{pa 3164 1559}%
\special{fp}%
%
\special{pn 20}%
\special{pa 3425 1602}%
\special{pa 3740 1293}%
\special{fp}%
%
\special{pn 20}%
\special{pa 3425 728}%
\special{pa 3704 994}%
\special{fp}%
%
\special{pn 20}%
\special{pa 1291 747}%
\special{pa 895 1320}%
\special{fp}%
%
\special{pn 20}%
\special{pa 1003 1484}%
\special{pa 1813 1484}%
\special{fp}%
%
\special{pn 20}%
\special{pa 1560 753}%
\special{pa 1875 1353}%
\special{fp}%
\put(6.8900,-15.6700){\makebox(0,0)[lb]{$x$}}%
\put(18.6800,-15.6700){\makebox(0,0)[lb]{$y$}}%
\put(25.7900,-12.0000){\makebox(0,0)[lb]{$y$}}%
\put(31.5500,-18.0800){\makebox(0,0)[lb]{$x$}}%
\put(52.3400,-12.0700){\makebox(0,0)[lb]{$y$}}%
\put(47.2100,-17.4800){\makebox(0,0)[lb]{$x$}}%
\put(8.6000,-4.4600){\makebox(0,0)[lb]{$A_2^{(1)}$}}%
\put(26.9000,-4.5300){\makebox(0,0)[lb]{$A_3^{(1)}$}}%
\put(48.8000,-4.3000){\makebox(0,0)[lb]{$D_4^{(1)}$}}%
\put(10.0000,-21.2500){\makebox(0,0)[lb]{$P_{IV}$}}%
\put(29.3000,-21.2500){\makebox(0,0)[lb]{$P_{V}$}}%
\put(50.4000,-21.1700){\makebox(0,0)[lb]{$P_{VI}$}}%
\end{picture}%
\label{fig:CPVI16}
\caption{Dynkin diagrams of types $A_2^{(1)},A_3^{(1)}$ and $D_4^{(1)}$}
\end{figure}

From the viewpoint of holomorphy, these transformations $u_1,u_2$ correspond to canonical coordinate systems $(x_i,y_i)$ $(i=1,2)$ (see \cite{Sasa5}), which are explicitly written as follows:
\begin{align*}
\begin{split}
&(x_1,y_1)=(1/x,-(xy+\gamma_1)x), \\ &(x_2,y_2)=(-(xy-\gamma_2)y,1/y).
\end{split}
\end{align*}
These canonical coordinate systems can be obtained by successive blowing-up procedures at the beginning of the accessible singular points
$$
P_1=\{(X_1,Y_1)=(0,0)\}, \ \ P_2=\{(X_2,Y_2)=(0,0)\}
$$
on the boundary divisor of ${\Bbb P}^2$. Here the coordinate systems $(X_i,Y_i)$ are the boundary coordinate systems of ${\Bbb P}^2$ with the rational transformations
$$
(X_1,Y_1)=(1/x,y/x), \qquad (X_2,Y_2)=(x/y,1/y).
$$

\begin{proposition}\label{P5.1}
Let us consider a polynomial Hamiltonian system with Hamiltonian $K \in {\Bbb C}(t)[x,y]$. We assume that

$(A1)$ $deg(K)=5$ with respect to $x,y$.

$(A2)$ This system becomes again a polynomial Hamiltonian system in each coordinate systems $(x_i,y_i) \ (i=1,2)$.

\noindent
Then such a system is explicitly given as follows{\rm : \rm}
\begin{equation}
  \left\{
  \begin{aligned}
   \frac{dx}{dt} &=\frac{\partial K}{\partial y}=2a_1x^3y+3a_2x^2y^2+2a_3x^2y+a_4x^2+2a_5xy+a_6x-\gamma_2a_5-\gamma_2^2 a_2,\\
   \frac{dy}{dt} &=-\frac{\partial K}{\partial x}=-3a_1x^2y^2-2a_2xy^3-2a_3xy^2-a_5y^2-2a_4xy-a_6y-\gamma_1a_4+{\gamma_1}^2a_1
   \end{aligned}
  \right. 
\end{equation}
with the polynomial Hamiltonian K
\begin{align}
\begin{split}
K=&a_1x^3y^2+a_2x^2y^3+a_3x^2y^2+a_4x^2y+a_5xy^2\\
&+a_6xy-(\gamma_2a_5+\gamma_2^2 a_2)y+(\gamma_1a_4-{\gamma_1}^2a_1)x.
\end{split}
\end{align}
Here, $a_1,a_2,\dots ,a_6$ are undetermined rational functions in t.
\end{proposition}

In the case of dimension four, it is easy to see that the affine Weyl groups $W(A_5^{(1)})$ and $W(A_4^{(1)})$ have a common subgroup $W$, whose elements $g_i$ are explicitly written as follows.
\begin{figure}[h]
\unitlength 0.1in
\begin{picture}(54.05,19.11)(6.00,-21.87)
%
\special{pn 20}%
\special{ar 1398 640 171 170  0.0000000 6.2831853}%
\put(8.4700,-4.4600){\makebox(0,0)[lb]{$A_4^{(1)}$}}%
\put(23.7600,-4.5000){\makebox(0,0)[lb]{$A_5^{(1)}$}}%
\put(44.3000,-4.5000){\makebox(0,0)[lb]{$D_6^{(1)}$}}%
%
\special{pn 20}%
\special{ar 771 1100 171 170  0.0000000 6.2831853}%
%
\special{pn 20}%
\special{ar 1151 1740 171 170  0.0000000 6.2831853}%
%
\special{pn 20}%
\special{ar 1797 1750 171 170  0.0000000 6.2831853}%
%
\special{pn 20}%
\special{ar 1996 1100 172 170  0.0000000 6.2831853}%
%
\special{pn 20}%
\special{ar 3012 687 171 170  0.0000000 6.2831853}%
%
\special{pn 20}%
\special{ar 2471 1117 171 170  0.0000000 6.2831853}%
%
\special{pn 20}%
\special{ar 2480 1647 171 170  0.0000000 6.2831853}%
%
\special{pn 20}%
\special{ar 3012 2017 171 170  0.0000000 6.2831853}%
%
\special{pn 20}%
\special{ar 3525 1117 171 170  0.0000000 6.2831853}%
%
\special{pn 20}%
\special{ar 3534 1647 172 170  0.0000000 6.2831853}%
%
\special{pn 20}%
\special{ar 3981 860 171 170  0.0000000 6.2831853}%
%
\special{pn 20}%
\special{ar 4418 1340 171 170  0.0000000 6.2831853}%
%
\special{pn 20}%
\special{ar 4019 1970 171 170  0.0000000 6.2831853}%
%
\special{pn 20}%
\special{ar 4950 1330 171 170  0.0000000 6.2831853}%
%
\special{pn 20}%
\special{ar 5425 1340 171 170  0.0000000 6.2831853}%
%
\special{pn 20}%
\special{ar 5795 880 172 170  0.0000000 6.2831853}%
%
\special{pn 20}%
\special{ar 5833 2000 172 170  0.0000000 6.2831853}%
%
\special{pn 20}%
\special{pa 1240 720}%
\special{pa 880 970}%
\special{fp}%
%
\special{pn 20}%
\special{pa 830 1250}%
\special{pa 1040 1590}%
\special{fp}%
%
\special{pn 20}%
\special{pa 1320 1750}%
\special{pa 1630 1750}%
\special{fp}%
%
\special{pn 20}%
\special{pa 1550 740}%
\special{pa 1850 970}%
\special{fp}%
%
\special{pn 20}%
\special{pa 2030 1280}%
\special{pa 1840 1580}%
\special{fp}%
%
\special{pn 20}%
\special{pa 2860 777}%
\special{pa 2570 987}%
\special{fp}%
%
\special{pn 20}%
\special{pa 2480 1287}%
\special{pa 2480 1477}%
\special{fp}%
%
\special{pn 20}%
\special{pa 2570 1787}%
\special{pa 2850 1927}%
\special{fp}%
%
\special{pn 20}%
\special{pa 3140 787}%
\special{pa 3420 977}%
\special{fp}%
%
\special{pn 20}%
\special{pa 3520 1277}%
\special{pa 3520 1467}%
\special{fp}%
%
\special{pn 20}%
\special{pa 3190 1967}%
\special{pa 3480 1817}%
\special{fp}%
%
\special{pn 20}%
\special{pa 4099 980}%
\special{pa 4309 1190}%
\special{fp}%
%
\special{pn 20}%
\special{pa 4139 1840}%
\special{pa 4389 1510}%
\special{fp}%
%
\special{pn 20}%
\special{pa 4599 1360}%
\special{pa 4779 1360}%
\special{fp}%
%
\special{pn 20}%
\special{pa 5119 1360}%
\special{pa 5249 1360}%
\special{fp}%
%
\special{pn 20}%
\special{pa 5539 1210}%
\special{pa 5719 1030}%
\special{fp}%
%
\special{pn 20}%
\special{pa 5499 1500}%
\special{pa 5709 1850}%
\special{fp}%
\put(6.5000,-11.9000){\makebox(0,0)[lb]{$y$}}%
\put(23.6000,-12.1700){\makebox(0,0)[lb]{$y$}}%
\put(42.8900,-14.3000){\makebox(0,0)[lb]{$y$}}%
\put(9.2000,-20.4000){\makebox(0,0)[lb]{$x-z$}}%
\put(22.4000,-19.1700){\makebox(0,0)[lb]{$x-z$}}%
\put(47.3900,-16.2000){\makebox(0,0)[lb]{$x-z$}}%
\put(16.7000,-18.4000){\makebox(0,0)[lb]{$w$}}%
\put(28.7000,-21.0700){\makebox(0,0)[lb]{$w$}}%
\put(52.9900,-14.5000){\makebox(0,0)[lb]{$w$}}%
\put(19.0000,-11.8000){\makebox(0,0)[lb]{$z$}}%
\put(34.1000,-17.4700){\makebox(0,0)[lb]{$z$}}%
\put(56.7900,-9.6000){\makebox(0,0)[lb]{$z$}}%
\end{picture}%
\label{fig:CPVI17}
\caption{Dynkin diagrams of types $A_4^{(1)},A_5^{(1)}$ and $D_6^{(1)}$}
\end{figure}
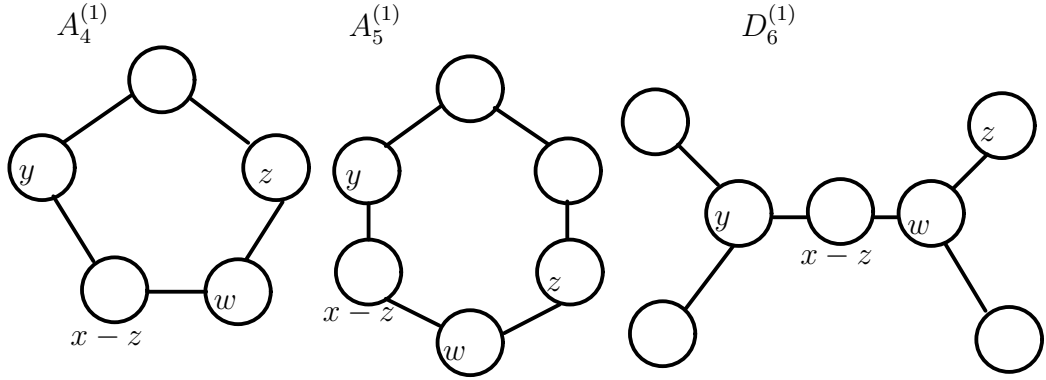
\begin{align*}
g_1: (x,y,z,w,t,\alpha_1,\ldots,\alpha_4) &\rightarrow \left(x,y,z+\frac{\alpha_1}{w},w,t,-\alpha_1,\alpha_2+\alpha_1,\alpha_3,\alpha_4+\alpha_1 \right), \\
g_2: (x,y,z,w,t,\alpha_1,\ldots,\alpha_4) &\rightarrow \left(x,y,z,w-\frac{\alpha_2}{z},t,\alpha_1+\alpha_2,-\alpha_2,\alpha_3,\alpha_4 \right), \\
g_3: (x,y,z,w,t,\alpha_1,\ldots,\alpha_4) &\rightarrow \left(x+\frac{\alpha_3}{y},y,z,w,t,\alpha_1,\alpha_2,-\alpha_3,\alpha_4+\alpha_3 \right), \\
g_4: (x,y,z,w,t,\alpha_1,\ldots,\alpha_4) &\rightarrow \left(x,y-\frac{\alpha_4}{x-z},z,w+\frac{\alpha_4}{x-z},t;\alpha_1+\alpha_4,\alpha_2,\alpha_3+\alpha_4,-\alpha_4 \right).
\end{align*}
Here, $\alpha_1,\alpha_2,\alpha_3$ and $\alpha_4$ are root parameters. 

\begin{proposition}\label{P5.2}
The transformations $g_i$ described above define a representation of the classical Weyl group of type $A_4$.
\end{proposition}

\begin{proposition}\label{P5.3}
Let us consider a polynomial Hamiltonian system with Hamiltonian $H \in {\Bbb C}(t)[x,y,z,w]$. We assume that

$(A1)$ $deg(H)=5$ with respect to $x,y,z,w$.

$(A2)$ This system becomes again a polynomial Hamiltonian system in each coordinate systems $(x_i,y_i,z_i,w_i)$ $(i=1,2,3,4)${\rm : \rm}
\begin{align*}
g_1:&x_1=x, \quad y_1=y, \quad z_1=1/z, \quad w_1=-z(zw+\alpha_1),\\
g_2:&x_2=x, \quad y_2=y, \quad z_2=-w(zw-\alpha_2), \quad w_2=1/w,\\
g_3:&x_3=1/x, \quad y_3=-x(xy+\alpha_3), \quad z_3=z, \quad w_3=w,\\
g_4:&x_4=-((x-z)y-\alpha_4)y, \quad y_4=1/y, \quad z_4=z, \quad w_4=y+w.
\end{align*}
Then such a system is explicitly given as follows{\rm : \rm}
\begin{equation}\label{37}
\frac{dx}{dt}=\frac{\partial H}{\partial y}, \ \ \frac{dy}{dt}=-\frac{\partial H}{\partial x}, \ \ \frac{dz}{dt}=\frac{\partial H}{\partial w}, \ \ \frac{dw}{dt}=-\frac{\partial H}{\partial z}
\end{equation}
with the polynomial Hamiltonian
\begin{equation}\label{38}
\begin{aligned}
&H=\frac{(b_5+b_9)}{2}x^3y^2+\frac{(b_4-b_7)}{2}x^2y^3+\frac{(b_1+b_3)}{2}x^2y^2+b_8x^2y+b_6xy^2+b_2xy\\
  &+\frac{(3\alpha_1\alpha_2+\alpha_2^2+2\alpha_2\alpha_4+\alpha_4^2)b_4+2(-\alpha_2-\alpha_4)b_6+(-\alpha_1\alpha_2-\alpha_2^2-2\alpha_2\alpha_4-\alpha_4^2)b_7}{2}y\\
  &-\frac{(\alpha_3b_5-2b_8+\alpha_3b_9)\alpha_3}{2}x+\frac{(b_5+b_9)}{2}z^3w^2+\frac{b_7-b_4}{2}z^2w^3+\frac{(b_1+b_3)}{2}z^2w^2\\
  &+\frac{(-3\alpha_1-4\alpha_4)b_4+2b_6+(\alpha_1+2\alpha_4)b_7}{2}zw^2+((\alpha_1+\alpha_4)b_1+b_2)zw\\
  &+\frac{(2\alpha_1+\alpha_4)b_5+2b_8+(-\alpha_4-2\alpha_3)b_9}{2}z^2w\\
  &+\frac{\alpha_1(\alpha_1+\alpha_4)b_5+2\alpha_1b_8+\alpha_1(-\alpha_1-2\alpha_3-\alpha_4)b_9}{2}z\\
  &+\frac{\alpha_2(3\alpha_1+\alpha_2+4\alpha_4)b_4-2\alpha_2b_6+\alpha_2(-\alpha_1-\alpha_2-2\alpha_4)b_7}{2}w+b_9(x^2yzw+\alpha_3xzw)\\
  &+b_1(yz^2w+\alpha_1yz)+b_3xyzw+b_5(xyz^2w+\alpha_1xyz)+2b_6yzw\\
  &+b_4(xyzw^2-\frac{5}{2}yz^2w^2-\frac{3}{2}y^2z^2w-3\alpha_1yzw-\frac{3}{2}\alpha_1y^2z-\alpha_2(x-z)yw-2\alpha_4yzw)\\
  &+b_7(xy^2zw+\frac{3}{2}yz^2w^2+\frac{1}{2}y^2z^2w+\alpha_1yzw+\frac{1}{2}\alpha_1y^2z+\alpha_4yzw).\\
\end{aligned}
\end{equation}
Here, $b_1,b_2,\dots ,,b_9$ are undetermined rational functions in t.
\end{proposition}
Propositions \ref{P5.1},\ref{P5.2} and \ref{P5.3} can be cheched by a direct calculation, respectively.

\section{Proof of Theorems \ref{1.3} and \ref{1.4}}

As is well-known, the degeneration from $P_{VI}$ to $P_{V}$ is given by
\begin{gather*}
\begin{gathered}
\alpha_0={\varepsilon}^{-1}, \ \alpha_1=A_3, \ \alpha_3=A_0-A_2-{\varepsilon}^{-1}, \ \alpha_4=A_1
\end{gathered}\\
t=1+{\varepsilon}T, \ (x-1)(X-1)=1, \ (x-1)y+(X-1)Y=-A_2.
\end{gather*}
Notice that
$$
A_0+A_1+A_2+A_3=\alpha_0+\alpha_1+2\alpha_2+\alpha_3+\alpha_4=1
$$
and the change of variables from $(q,p)$ to $(Q,P)$ is symplectic.

As the fourth-order analogue of the above confluence process, we consider the following coupling confluence process from the system \eqref{1}. We take the following coupling confluence process $P_{VI} \rightarrow P_{V}$ for each coordinate system $(x,y)$ and $(z,w)$ of the system \eqref{1}
\begin{gather*}
\begin{gathered}
\alpha_0={\varepsilon}^{-1}, \quad \alpha_1=A_0, \quad \alpha_2=A_1, \quad \alpha_4-\beta_4=A_2,\\
\beta_2=A_3, \quad \beta_3=-{\varepsilon}^{-1}-(A_1+A_2+A_3-A_5), \quad \beta_4=A_4,
\end{gathered}\\
\begin{gathered}
t=1-{\varepsilon}T, \quad x=\frac{X}{X-T}, \quad z=\frac{Z}{Z-T},\\
y=-\frac{(X-T)\{(X-T)Y+A_1\}}{T}, \quad w=-\frac{(Z-T)\{(Z-T)W+A_3\}}{T}
\end{gathered}
\end{gather*}
from $\alpha_0,\alpha_1,\alpha_2,\gamma_1,\beta_2,\beta_3,\beta_4,t,x,y,z,w$ to $A_0,\ldots ,A_5,\varepsilon,T,X,Y,Z,W$. Notice that
$$
A_0+A_1+A_2+A_3+A_4+A_5=\alpha_0+\alpha_1+2\alpha_2+2(\alpha_4-\beta_4)+2\beta_2+\beta_3+\beta_4=1
$$
and the change of variables from $(x,y,z,w)$ to $(X,Y,Z,W)$ is symplectic. Then the system \eqref{1} can also be written in the new variables $T,X,Y,Z,W$ and parameters $A_0,A_1,\dots ,A_5,\varepsilon$ as a Hamiltonian system. This new system tends to the system \eqref{29} of type $A_5^{(1)}$ as $\varepsilon \rightarrow 0$. The proof of Theorem \ref{1.3} is now complete.

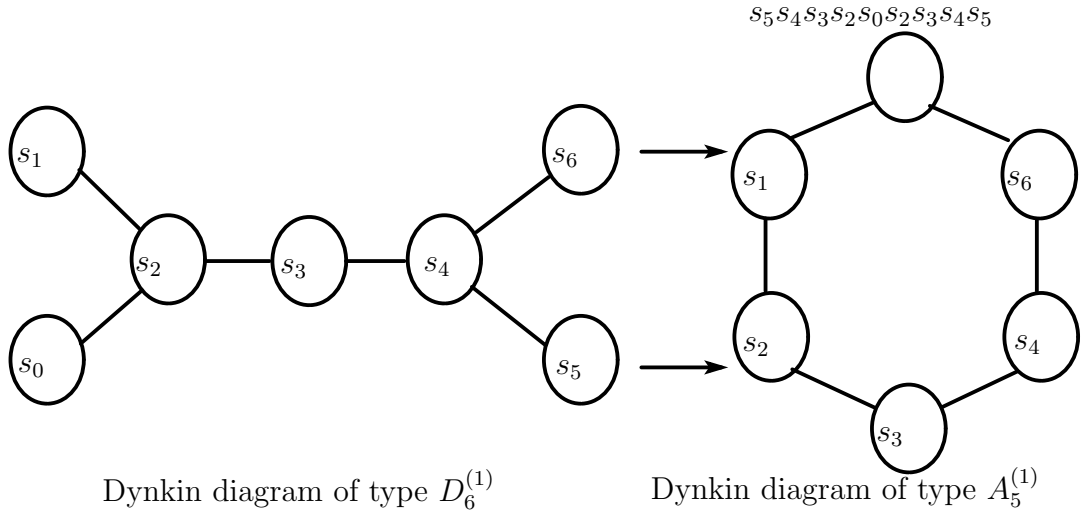
\begin{figure}[h]
\unitlength 0.1in
\begin{picture}(58.86,25.30)(9.40,-36.20)
%
\special{pn 20}%
\special{ar 1430 1919 192 230  1.5707963 6.2831853}%
\special{ar 1430 1919 192 230  0.0000000 1.5291537}%
%
\special{pn 20}%
\special{ar 1430 3010 192 230  1.5707963 6.2831853}%
\special{ar 1430 3010 192 230  0.0000000 1.5291537}%
%
\special{pn 20}%
\special{ar 2064 2484 192 231  1.5707963 6.2831853}%
\special{ar 2064 2484 192 231  0.0000000 1.5291537}%
%
\special{pn 20}%
\special{ar 3509 2484 192 231  1.5707963 6.2831853}%
\special{ar 3509 2484 192 231  0.0000000 1.5291537}%
%
\special{pn 20}%
\special{ar 4223 1910 193 230  1.5707963 6.2831853}%
\special{ar 4223 1910 193 230  0.0000000 1.5293693}%
%
\special{pn 20}%
\special{ar 4223 3010 193 230  1.5707963 6.2831853}%
\special{ar 4223 3010 193 230  0.0000000 1.5293693}%
%
\special{pn 20}%
\special{pa 1598 2016}%
\special{pa 1911 2322}%
\special{fp}%
%
\special{pn 20}%
\special{pa 1606 2925}%
\special{pa 1927 2646}%
\special{fp}%
%
\special{pn 20}%
\special{ar 2802 2493 193 231  1.5707963 6.2831853}%
\special{ar 2802 2493 193 231  0.0000000 1.5293693}%
%
\special{pn 20}%
\special{pa 2256 2493}%
\special{pa 2594 2493}%
\special{fp}%
%
\special{pn 20}%
\special{pa 2995 2493}%
\special{pa 3300 2493}%
\special{fp}%
%
\special{pn 20}%
\special{pa 3677 2341}%
\special{pa 4055 2053}%
\special{fp}%
%
\special{pn 20}%
\special{pa 3653 2628}%
\special{pa 4039 2934}%
\special{fp}%
%
\put(9.4000,-19.8600){\makebox(0,0)[lb]{}}%
\put(12.7700,-30.9800){\makebox(0,0)[lb]{$s_0$}}%
\put(18.9000,-25.7000){\makebox(0,0)[lb]{$s_2$}}%
\put(26.5000,-25.8000){\makebox(0,0)[lb]{$s_3$}}%
\put(34.0000,-25.7000){\makebox(0,0)[lb]{$s_4$}}%
\put(40.9000,-31.1000){\makebox(0,0)[lb]{$s_5$}}%
\put(40.7000,-20.0000){\makebox(0,0)[lb]{$s_6$}}%
\put(12.7000,-20.0000){\makebox(0,0)[lb]{$s_1$}}%
%
\special{pn 20}%
\special{pa 4540 1920}%
\special{pa 4960 1920}%
\special{fp}%
\special{sh 1}%
\special{pa 4960 1920}%
\special{pa 4893 1900}%
\special{pa 4907 1920}%
\special{pa 4893 1940}%
\special{pa 4960 1920}%
\special{fp}%
%
\special{pn 20}%
\special{pa 4540 3050}%
\special{pa 4960 3050}%
\special{fp}%
\special{sh 1}%
\special{pa 4960 3050}%
\special{pa 4893 3030}%
\special{pa 4907 3050}%
\special{pa 4893 3070}%
\special{pa 4960 3050}%
\special{fp}%
%
\special{pn 20}%
\special{ar 5920 1530 193 230  1.5707963 6.2831853}%
\special{ar 5920 1530 193 230  0.0000000 1.5293693}%
%
\special{pn 20}%
\special{ar 5210 2040 193 230  1.5707963 6.2831853}%
\special{ar 5210 2040 193 230  0.0000000 1.5293693}%
%
\special{pn 20}%
\special{ar 5220 2890 193 230  1.5707963 6.2831853}%
\special{ar 5220 2890 193 230  0.0000000 1.5293693}%
%
\special{pn 20}%
\special{ar 5940 3370 193 230  1.5707963 6.2831853}%
\special{ar 5940 3370 193 230  0.0000000 1.5293693}%
%
\special{pn 20}%
\special{ar 6623 2040 193 230  1.5707963 6.2831853}%
\special{ar 6623 2040 193 230  0.0000000 1.5293693}%
%
\special{pn 20}%
\special{ar 6633 2890 193 230  1.5707963 6.2831853}%
\special{ar 6633 2890 193 230  0.0000000 1.5293693}%
%
\special{pn 20}%
\special{pa 5760 1660}%
\special{pa 5320 1850}%
\special{fp}%
%
\special{pn 20}%
\special{pa 5190 2270}%
\special{pa 5190 2650}%
\special{fp}%
%
\special{pn 20}%
\special{pa 5330 3060}%
\special{pa 5760 3260}%
\special{fp}%
%
\special{pn 20}%
\special{pa 6050 1680}%
\special{pa 6460 1860}%
\special{fp}%
%
\special{pn 20}%
\special{pa 6610 2280}%
\special{pa 6610 2650}%
\special{fp}%
%
\special{pn 20}%
\special{pa 6110 3260}%
\special{pa 6510 3070}%
\special{fp}%
\put(50.7000,-21.2000){\makebox(0,0)[lb]{$s_1$}}%
\put(50.7000,-29.9000){\makebox(0,0)[lb]{$s_2$}}%
\put(57.7000,-34.6000){\makebox(0,0)[lb]{$s_3$}}%
\put(64.9000,-29.8000){\makebox(0,0)[lb]{$s_4$}}%
\put(64.5000,-21.2000){\makebox(0,0)[lb]{$s_6$}}%
\put(51.0000,-12.6000){\makebox(0,0)[lb]{$s_5s_4s_3s_2s_0s_2s_3s_4s_5$}}%
\put(17.2000,-37.9000){\makebox(0,0)[lb]{Dynkin diagram of type $D_6^{(1)}$}}%
\put(45.9000,-37.8000){\makebox(0,0)[lb]{Dynkin diagram of type $A_5^{(1)}$}}%
\end{picture}%
\label{fig:D64}
\caption{Confluence process from the system of $D_6^{(1)}$ to the system of $A_5^{(1)}$}
\end{figure}

Next, let us prove Theorem \ref{1.4}. Choose $S_i$, $i=0,1,2,3,4,5$ as
$$
S_0:=s_1, \ S_1:=s_2, \ S_2:=s_3, \ S_3:=s_4, \ S_4:=s_6, \ S_5:=s_5s_4s_3s_2s_0s_2s_3s_4s_5,
$$
which are reflections of
\begin{gather*}
\begin{gathered}
A_0:=\alpha_1, \quad A_1:=\alpha_2, \quad A_2:=\alpha_4-\beta_4, \quad A_3:=\beta_2,\\
A_4:=\beta_4, \quad A_5:=1-\alpha_1-\alpha_2-\alpha_4-\beta_2.
\end{gathered}
\end{gather*}
By using the notation $(*):=(A_0,A_1,A_2,A_3,A_4,A_5,\varepsilon)$, we can easily check
\begin{align*}
S_0(*)&=(-A_0,A_1+A_0,A_2,A_3,A_4,A_5+A_0,\varepsilon),\\
S_1(*)&=\left(A_0+A_1,-A_1,A_2+A_1,A_3,A_4,A_5,\frac{\varepsilon}{1+\varepsilon A_1} \right),\\
S_2(*)&=(A_0,A_1+A_2,-A_2,A_3+A_2,A_4,A_5,\varepsilon),\\
S_3(*)&=(A_0,A_1,A_2+A_3,-A_3,A_4+A_3,A_5,\varepsilon),\\
S_4(*)&=(A_0,A_1,A_2,A_3+A_4,-A_4,A_5+A_4,\varepsilon),\\
S_5(*)&=\left(A_0+A_5,A_1,A_2,A_3,A_4+A_5,-A_5,\frac{\varepsilon}{1-\varepsilon A_5} \right).
\end{align*}
By the above relation, we will see that the group $<S_0,S_1,\dots,S_5>$ can be considered to be an affine Weyl group of the affine Lie algebra of type $A_5^{(1)}$ with respect to simple roots $A_0,A_1,\dots,A_5$.

Now we investigate how the generators of $<S_0,S_1,\dots,S_5>$ act on $X,Y,Z,W$. By using the notation $(**):=(X,Y,Z,W)$, we can verify
\begin{align*}
S_0(**)&=\left(X,Y-\frac{A_0}{X-T},Z,W \right),\\
S_1(**)&=\left(X+\frac{A_1}{Y},Y,Z,W \right),\\
S_2(**)&=\left(X,Y-\frac{A_2}{X-Z},Z,W+\frac{A_2}{X-Z} \right),\\
S_3(**)&=\left(X,Y,Z+\frac{A_3}{W},W \right),\\
S_4(**)&=\left(X,Y,Z,W-\frac{A_4}{Z} \right),\\
S_5(**)&=\left(X+\frac{A_5}{Y+W-1},Y,Z+\frac{A_5}{Y+W-1},W \right).
\end{align*}
The proof of Theorem \ref{1.4} has thus been completed. \qed

\section{Appendix }
For the second-order Painlev\'e equations, we can obtain the entire space of initial conditions by adding subvarieties of codimension 1 (equivalently, of dimension 1) to the space of initial conditions of holomorphic solutions. However, in the case of the fourth-order differential equations, we need to add codimension 2 subvarieties to the space in addition to codimension 1 subvarieties (see \cite{Sasa10,Sasa11,Sasa12}). In order to resolve singularities, we need to both blow up and blow down. Moreover, to obtain a smooth variety by blowing-down, we need to resolve for a pair of singularities (see \cite{Sasa10,Sasa11,Sasa12}). In this section, we will give some canonical coordinate systems of the system \eqref{1}. Each of them contains a 3-parameter or 2-parameter family of meromorphic solutions.

In order to consider the singularity analysis for the system \eqref{1}, let us take the compactification $([z_0:z_1:z_2:z_3:z_4],t) \in {\Bbb P}^4 \times B$ of $(x,y,z,w,t) \in {\Bbb C}^4 \times B$, $B:={\Bbb C}-\{0,1\}$ with the natural embedding
$$
(x,y,z,w)=(z_1/z_0,z_2/z_0,z_3/z_0,z_4/z_0).
$$
Moreover, we denote the boundary divisor in ${\Bbb P}^4$ by $ {\mathcal H}$. Fixing the parameter ${\alpha}_0,\alpha_1,\beta_0$, consider the product ${\Bbb P}^4 \times B$ and extend the regular vector field on ${\Bbb C}^4 \times B$ to a rational vector field $\tilde v$ on ${\Bbb P}^4 \times B$. It is easy to see that ${\Bbb P}^4 \times B$ is covered by five copies of ${\Bbb C}^4 \times B${\rm : \rm}
\begin{gather*}
\begin{gathered}
U_0 \times B={\Bbb C}^4 \times B \ni (x,y,z,w,t),
\end{gathered}\\
U_j \times B={\Bbb C}^4 \times B \ni (X_j,Y_j,Z_j,W_j,t) \ (j=1,2,3,4)
\end{gather*}
via the following rational transformations
\begin{align*}
1)&X_1=1/x, \quad Y_1=y/x, \quad Z_1=z/x, \quad W_1=w/x,\\
2)&X_2=x/z, \quad Y_2=y/z, \quad Z_2=1/z, \quad W_2=w/z,\\
3)&X_3=x/y, \quad Y_3=1/y, \quad Z_3=z/y, \quad W_3=w/y,\\
4)&X_4=x/w, \quad Y_4=y/w, \quad Z_4=z/w, \quad W_4=1/w.
\end{align*}

The following lemma shows that this rational vector field $\tilde v$ has the following five accessible singular loci on the boundary divisor $\mathcal H \times \{ $t$ \} \subset {\Bbb P}^4 \times \{ $t$ \}$ for each $t \in B$.

\begin{lemma}\label{9.1}

The rational vector field $\tilde v$ has the following accessible singular loci{\rm : \rm}
\begin{equation*}
  \left\{
  \begin{aligned}
   P_i &=\{(X_i,Y_i,Z_i,W_i)|X_i=Y_i=Z_i=W_i=0\} \ (1=1,2,3,4),\\
   P_5 &=\{(X_3,Y_3,Z_3,W_3)|X_3=Y_3=Z_3=0,W_3=-1\}.\\
   \end{aligned}
  \right. 
\end{equation*}

\end{lemma}

This lemma can be proven by a direct calculation. \qed

\begin{proposition}\label{pro:ini}
If we resolve the accessible singular points given in Lemma \ref{9.1} by blowing-ups, then we can obtain the canonical coordinates $r_j (j=0,1,\dots ,6)$ given in Theorem \ref{1.2}.
\end{proposition}

{\it Proof} \hspace{0.2cm} By the following steps, we can resolve the accessible singular point $P_1$.

{\bf Step 1}: We blow up at the point $P_1${\rm : \rm}
$$
{X_1}^{(1)}=X_1 \;, \;\;\; {Y_1}^{(1)}=\frac{Y_1}{X_1} \;, \;\;\; {Z_1}^{(1)}=\frac{Z_1}{X_1} \;, \;\;\; {W_1}^{(1)}=\frac{W_1}{X_1}.
$$

{\bf Step 2}: We blow up along the surface $\{({X_1}^{(1)},{Y_1}^{(1)},{Z_1}^{(1)},{W_1}^{(1)})| {X_1}^{(1)}={Y_1}^{(1)}\\
=0\}${\rm : \rm}
$$
{X_1}^{(2)}={X_1}^{(1)} \;, \;\;\; {Y_1}^{(2)}=\frac{{Y_1}^{(1)}}{{X_1}^{(1)}} \;, \;\;\; {Z_1}^{(2)}={Z_1}^{(1)} \;, \;\;\; {W_1}^{(2)}={W_1}^{(1)}.
$$

It is easy to see that there are two accessible singular loci:
\begin{align*}
S_1^{(1)}&=\{({X_1}^{(2)},{Y_1}^{(2)},{Z_1}^{(2)},{W_1}^{(2)})|{X_1}^{(2)}={Y_1}^{(2)}+\alpha_1+\alpha_2=0 \},\\
S_1^{(2)}&=\{({X_1}^{(2)},{Y_1}^{(2)},{Z_1}^{(2)},{W_1}^{(2)})|{X_1}^{(2)}={Y_1}^{(2)}+\alpha_2=0 \}.
\end{align*}

{\bf Step 3}: We blow up along the surface $S_1^{(1)}${\rm : \rm}
$$
{X_1}^{(3)}={X_1}^{(2)} \;, \;\;\; {Y_1}^{(3)}=\frac{{Y_1}^{(2)}+\alpha_1+\alpha_2}{{X_1}^{(2)}} \;, \;\;\; {Z_1}^{(3)}={Z_1}^{(2)} \;, \;\;\; {W_1}^{(3)}={W_1}^{(2)}.
$$

{\bf Step 4}: We blow up along the surface $S_1^{(2)}${\rm : \rm}
$$
{X_1}^{(4)}={X_1}^{(2)} \;, \;\;\; {Y_1}^{(4)}=\frac{{Y_1}^{(2)}+\alpha_2}{{X_1}^{(2)}} \;, \;\;\; {Z_1}^{(4)}={Z_1}^{(2)} \;, \;\;\; {W_1}^{(4)}={W_1}^{(2)}.
$$
Thus we have resolved the accessible singular point $P_1$. 

By choosing new coordinate systems as
$$
(x_k,y_k,z_k,w_k)=({X_1}^{(k+2)},-{Y_1}^{(k+2)},{Z_1}^{(k+2)},{W_1}^{(k+2)}) \ (k=1,2),
$$
we can obtain the coordinate systems $(x_k,y_k,z_k,w_k) \ (k=1,2)$, respectively.

\vspace{0.5cm}
By the following steps, we can resolve the accessible singular point $P_4$.

{\bf Step 1}: We blow up at the point $P_4${\rm : \rm}
$$
{X_4}^{(1)}=\frac{X_4}{W_4} \;, \;\;\; {Y_4}^{(1)}=\frac{Y_4}{W_4} \;, \;\;\; {Z_4}^{(1)}=\frac{Z_4}{W_4} \;, \;\;\; {W_4}^{(1)}=W_4.
$$

It is easy to see that there are two accessible singular loci:
\begin{align*}
S_4^{(1)}&=\{({X_4}^{(1)},{Y_4}^{(1)},{Z_4}^{(1)},{W_4}^{(1)})|{Z_4}^{(1)}-1={W_4}^{(1)}=0 \},\\
S_4^{(2)}&=\{({X_4}^{(1)},{Y_4}^{(1)},{Z_4}^{(1)},{W_4}^{(1)})|{Z_4}^{(1)}={W_4}^{(1)}=0 \}.
\end{align*}

{\bf Step 2}: We blow up along the surface $S_4^{(1)}${\rm : \rm}
$$
{X_4}^{(2)}={X_4}^{(1)} \;, \;\;\; {Y_4}^{(2)}={Y_4}^{(1)} \;, \;\;\; {Z_4}^{(2)}=\frac{{Z_4}^{(1)}-1}{{W_4}^{(1)}} \;, \;\;\; {W_4}^{(2)}={W_4}^{(1)}.
$$

{\bf Step 3}: We blow up along the surface $\{({X_4}^{(2)},{Y_4}^{(2)},{Z_4}^{(2)},{W_4}^{(2)})|{Z_4}^{(2)}-\beta_3\\
={W_4}^{(2)}=0 \}${\rm : \rm}
$$
{X_4}^{(3)}={X_4}^{(2)} \;, \;\;\; {Y_4}^{(3)}={Y_4}^{(2)} \;, \;\;\; {Z_4}^{(3)}=\frac{{Z_4}^{(2)}-\beta_3}{{W_4}^{(2)}} \;, \;\;\; {W_4}^{(3)}={W_4}^{(2)}.
$$

{\bf Step 4}: We blow up along the surface $S_4^{(2)}${\rm : \rm}
$$
{X_4}^{(4)}={X_4}^{(1)} \;, \;\;\; {Y_4}^{(4)}={Y_4}^{(1)} \;, \;\;\; {Z_4}^{(4)}=\frac{{Z_4}^{(1)}}{{W_4}^{(1)}} \;, \;\;\; {W_4}^{(4)}={W_4}^{(1)}.
$$

{\bf Step 5}: We blow up along the surface $\{({X_4}^{(4)},{Y_4}^{(4)},{Z_4}^{(4)},{W_4}^{(4)})|{Z_4}^{(4)}-\beta_4\\
={W_4}^{(4)}=0 \}${\rm : \rm}
$$
{X_4}^{(5)}={X_4}^{(4)} \;, \;\;\; {Y_4}^{(5)}={Y_4}^{(4)} \;, \;\;\; {Z_4}^{(5)}=\frac{{Z_4}^{(4)}-\beta_4}{{W_4}^{(4)}} \;, \;\;\; {W_4}^{(5)}={W_4}^{(4)}.
$$
Thus we have resolved the accessible singular point $P_4$. 

By choosing new coordinate systems as
\begin{align*}
(x_5,y_5,z_5,w_5)&=({X_4}^{(3)},{Y_4}^{(3)},-{Z_4}^{(3)},{W_4}^{(3)})\\
(x_6,y_6,z_6,w_6)&=({X_4}^{(5)},{Y_4}^{(5)},-{Z_4}^{(5)},{W_4}^{(5)}),
\end{align*}
we can obtain the coordinate systems $(x_k,y_k,z_k,w_k) \ (k=5,6)$, respectively.

\vspace{0.5cm}
By the following steps, we can resolve the accessible singular point $P_5$.

{\bf Step 0}: We take the coordinate system centered at $P_5${\rm : \rm}
$$
{X_5}^{(0)}=X_3 \;, \;\;\; {Y_5}^{(0)}=Y_3 \;, \;\;\; {Z_5}^{(0)}=Z_3 \;, \;\;\; {W_5}^{(0)}=W_3+1.
$$

{\bf Step 1}: We blow up at the point $P_5${\rm : \rm}
$$
{X_5}^{(1)}=\frac{{X_5}^{(0)}}{{Y_5}^{(0)}} \;, \;\;\; {Y_5}^{(1)}={Y_5}^{(0)} \;, \;\;\; {Z_5}^{(1)}=\frac{{Z_5}^{(0)}}{{Y_5}^{(0)}} \;, \;\;\; {W_5}^{(1)}=\frac{{W_5}^{(0)}}{{Y_5}^{(0)}}.
$$

{\bf Step 2}: We blow up along the surface $\{({X_5}^{(1)},{Y_5}^{(1)},{Z_5}^{(1)},{W_5}^{(1)})| {X_5}^{(1)}-{Z_5}^{(1)}\\
={Y_5}^{(1)}=0\}${\rm : \rm}
$$
{X_5}^{(2)}=\frac{{X_5}^{(1)}-{Z_5}^{(1)}}{{Y_5}^{(1)}} \;, \;\;\; {Y_5}^{(2)}={Y_5}^{(1)} \;, \;\;\; {Z_5}^{(2)}={Z_5}^{(1)} \;, \;\;\; {W_5}^{(2)}={W_5}^{(1)}.
$$

{\bf Step 3}: We blow up along the surface $\{({X_5}^{(2)},{Y_5}^{(2)},{Z_5}^{(2)},{W_5}^{(2)})|{X_5}^{(2)}-\\
(\alpha_4-\beta_4)={Y_5}^{(2)}=0 \}${\rm : \rm}
$$
{X_5}^{(3)}=\frac{{X_5}^{(2)}-(\alpha_4-\beta_4)}{{Y_5}^{(2)}} \;, \;\;\; {Y_5}^{(3)}={Y_5}^{(2)} \;, \;\;\; {Z_5}^{(3)}={Z_5}^{(2)} \;, \;\;\; {W_5}^{(3)}={W_5}^{(2)}.
$$
Thus we have resolved the accessible singular point $P_5$. 

By choosing a new coordinate system as
$$
(x_3,y_3,z_3,w_3)=(-{X_5}^{(3)},{Y_5}^{(3)},{Z_5}^{(3)},{W_5}^{(3)}),
$$
we can obtain the coordinate system $(x_3,y_3,z_3,w_3)$.

For the remaining accessible singular points, the proof is similar.

Collecting all the cases, we have obtained the canonical coordinate systems $(x_j,y_j,\\
z_j,w_j) \ (j=0,1,\dots ,6)$, which proves Proposition 10.1.  \qed

We remark that each coordinate system contains a three-parameter family of meromorphic solutions of \eqref{1} as the initial conditions.

\vspace{0.2cm}
By using the coordinate systems $(x_j,y_j,z_j,w_j)$ $(j=0,1,2,\dots ,6)$, we will now make coordinate systems associated with other small meromorphic solution spaces. For example, we can take the coordinate system $(x_3,y_3,z_3,w_3)=(-((x-z)y-(\alpha_4-\beta_4))y,1/y,z,y+w)$. As a boundary coordinate system of this system $(x_3,y_3,z_3,w_3)$, we can take the coordinate system
$$
(X_3^{(0)},Y_3^{(0)},Z_3^{(0)},W_3^{(0)})=(x_3,y_3,z_3,1/w_3).
$$
It is easy to see that there is an accessible singularity along the surface
$$
S_3=\{(X_3^{(0)},Y_3^{(0)},Z_3^{(0)},W_3^{(0)})|Z_3^{(0)}=W_3^{(0)}=0\}.
$$
Now we blow up along the accessible singularity $S_3$.

{\bf Step 1}: We blow up along the surface $S_3${\rm : \rm}
$$
{X_3}^{(1)}=X_3^{(0)} \;, \;\;\; {Y_3}^{(1)}=Y_3^{(0)} \;, \;\;\; {Z_3}^{(1)}=\frac{Z_3^{(0)}}{W_3^{(0)}} \;, \;\;\; {W_3}^{(1)}=W_3^{(0)}.
$$

{\bf Step 2}: We blow up along the surface $\{(X_3^{(1)},Y_3^{(1)},Z_3^{(1)},W_3^{(1)})| Z_3^{(1)}-\beta_4=W_3^{(1)}=0\}${\rm : \rm}
$$
{X_3}^{(2)}={X_3}^{(1)} \;, \;\;\; {Y_3}^{(2)}={Y_3}^{(1)} \;, \;\;\; {Z_3}^{(2)}=\frac{{Z_3}^{(1)}-\beta_4}{{W_3}^{(1)}} \;, \;\;\; {W_3}^{(2)}={W_3}^{(1)}.
$$
Thus we have resolved the accessible singularity $S_3$. By the same way, we can obtain the following canonical coordinate systems.

\begin{proposition}
The system \eqref{1} has the following canonical coordinate systems with regard to the transformations $r_ir_j${\rm : \rm}

\begin{align*}
\begin{split}
r_0r_3:&x_7=-(y+w)((x-t)(y+w)-\alpha_0), \quad y_7=1/(y+w),\\
&z_7=-((z-x)w-(\alpha_4-\beta_4))w, \quad w_7=1/w,
\end{split}\\
\begin{split}
r_0r_3:&x_7=-(y+w)((x-t)(y+w)-\alpha_0), \quad y_7=1/(y+w),\\
&z_7=-((z-x)w-(\alpha_4-\beta_4))w, \quad w_7=1/w,
\end{split}\\
\begin{split}
r_0r_4:&x_{8}=-((x-t)y-\alpha_0)y, \quad y_{8}=1/y,\\
&z_{8}=1/z, \quad w_{8}=-z(zw+\beta_2),
\end{split}\\
\begin{split}
r_0r_5:&x_9=-((x-t)y-\alpha_0)y, \quad y_9=1/y,\\
&z_9=-((z-1)w-\beta_3)w, \quad w_9=1/w,
\end{split}\\
\begin{split}
r_0r_6:&x_{10}=-((x-t)y-\alpha_0)y, \quad y_{10}=1/y,\\
&z_{10}=-(zw-\beta_4)w, \quad w_{10}=1/w,
\end{split}\\
\begin{split}
r_1r_4:&x_{11}=1/x, \quad y_{11}=-(xy+\alpha_1+\alpha_2)x,\\
&z_{11}=1/z, \quad w_{11}=-z(zw+\beta_2),
\end{split}\\
\begin{split}
r_1r_5:&x_{12}=1/x, \quad y_{12}=-(xy+\alpha_1+\alpha_2)x,\\
&z_{12}=-((z-1)w-\beta_3)w, \quad w_{12}=1/w,
\end{split}
\end{align*}

\begin{align*}
\begin{split}
r_1r_6:&x_{13}=1/x, \quad y_{13}=-(xy+\alpha_1+\alpha_2)x,\\
&z_{13}=-(zw-\beta_4)w, \quad w_{13}=1/w,
\end{split}\\
\begin{split}
r_2r_4:&x_{14}=1/x, \quad y_{14}=-x(xy+\alpha_2),\\
&z_{14}=1/z, \quad w_{14}=-z(zw+\beta_2),
\end{split}\\
\begin{split}
r_2r_5:&x_{15}=1/x, \quad y_{15}=-x(xy+\alpha_2),\\
&z_{15}=-((z-1)w-\beta_3)w, \quad w_{15}=1/w,
\end{split}\\
\begin{split}
r_2r_6:&x_{16}=1/x, \quad y_{16}=-x(xy+\alpha_2),\\
&z_{16}=-(zw-\beta_4)w, \quad w_{16}=1/w,
\end{split}\\
\begin{split}
r_3r_5:&x_{17}=-((x-z)y-(\alpha_4-\beta_4))y, \quad y_{17}=1/y,\\
&z_{17}=-((z-1)(y+w)-\beta_3)(y+w), \quad w_{17}=1/(y+w),
\end{split}\\
\begin{split}
r_3r_6:&x_{18}=-((x-z)y-(\alpha_4-\beta_4))y, \quad y_{18}=1/y,\\
&z_{18}=-(z(y+w)-\beta_4)(y+w), \quad w_{18}=1/(y+w).
\end{split}
\end{align*}
\end{proposition}
Each coordinate system contains a two-parameter family of meromorphic solutions of \eqref{1} as the initial conditions. By using the coordinate systems $(x_j,y_j,z_j,w_j)$ $(j=7,8,\dots ,18)$, we will now make the coordinate systems associated with other small meromorphic solution spaces by the same way. For example, we can take the coordinate system $(x_{15},y_{15},z_{15},w_{15})=(-((x-z)y-(\alpha_4-\beta_4))y,1/y,-(z(y+w)-\beta_4)(y+w),1/(y+w))$. As a boundary coordinate system of this system $(x_{15},y_{15},z_{15},w_{15})$, we can take the coordinate system
$$
(X_{15},Y_{15},Z_{15},W_{15})=(x_{15}+z_{15},y_{15},1/z_{15},w_{15}-y_{15}).
$$
It is easy to see that there is an accessible singularity along the surface
$$
S_{15}=\{(X_{15},Y_{15},Z_{15},W_{15})|Z_{15}=W_{15}=0\}.
$$
Now we blow up along the accessible singularity $S_{15}$.

{\bf Step 1}: We blow up along the surface $\{(X_{15},Y_{15},Z_{15},W_{15})| Z_{15}=W_{15}=0\}${\rm : \rm}
$$
{X_{15}}^{(1)}=X_{15} \;, \;\;\; {Y_{15}}^{(1)}=Y_{15} \;, \;\;\; {Z_{15}}^{(1)}=Z_{15} \;, \;\;\; {W_{15}}^{(1)}=\frac{W_{15}}{Z_{15}}.
$$

{\bf Step 2}: We blow up along the surface $\{({X_{15}}^{(1)},{Y_{15}}^{(1)},{Z_{15}}^{(1)},{W_{15}}^{(1)})|{Z_{15}}^{(1)}={W_{15}}^{(1)}+\beta_2=0 \}${\rm : \rm}
$$
{X_{15}}^{(2)}={X_{15}}^{(1)} \;, \;\;\; {Y_{15}}^{(2)}={Y_{15}}^{(1)} \;, \;\;\; {Z_{15}}^{(2)}={Z_{15}}^{(1)} \;, \;\;\; {W_{15}}^{(2)}=\frac{{W_{15}}^{(1)}+\beta_2}{{Z_{15}}^{(1)}}.
$$
Thus we have resolved the accessible singularity $S_{15}$. By the same way, we can obtain the following canonical coordinate systems.

\begin{proposition}
The system \eqref{1} has the following canonical coordinate systems with regard to the transformations $r_ir_jr_k${\rm : \rm}

\begin{align*}
\begin{split}
&r_3(r_4r_2):\\
&x_{19}=1/x, \quad y_{19}=-x^2y-z^2w-\alpha_2x-\beta_2z, \\
&z_{19}=z(zw+\beta_2)(-xzw+z^2w-\alpha_4x-\beta_2x+\beta_2z+\beta_4x)/x,\\
&w_{19}=-1/(z(zw+\beta_2)),
\end{split}\\
\begin{split}
&r_4(r_5r_3):\\
&x_{20}=w^2+2yw+y^2-xy^2-zw^2-2yzw+\alpha_4y+\beta_3(y+w)-\beta_4y, \\
&y_{20}=1/y, \quad z_{20}=-1/((y+w)(-y-w+zw+yz-\beta_3)), \\
&w_{20}=-((y+w)(-y-w+zw+yz-\beta_3)((y+w)(zw-w)+\beta_2y-\beta_3w))/y,
\end{split}\\
\begin{split}
&r_4(r_6r_3):\\
&x_{21}=-xy^2-zw^2-2yzw+\alpha_4y+\beta_4w, \quad y_{21}=1/y,\\
&z_{21}=-1/((y+w)(zw+yz-\beta_4)),\\
&w_{21}=-(y+w)(zw+yz-\beta_4)(zw^2+yzw+\beta_2y-\beta_4w)/y,
\end{split}\\
\begin{split}
&r_2(r_0r_3):\\
&x_{22}=1/((y+w)(tw-xw+ty-xy+\alpha_0)),\\
&y_{22}=(y+w)(tw-xw+ty-xy+\alpha_0)(tyw-xyw+ty^2-xy^2+\alpha_0y-\alpha_2w)/w,\\
&z_{22}=tw^2+2(t-x)yw+ty^2-xy^2-zw^2+\alpha_0(y+w)+(\alpha_4-\beta_4)w, \ w_{22}=1/w,
\end{split}\\
\begin{split}
&r_3(r_4r_1):\\
&x_{23}=1/x, \quad y_{23}=-x^2y-z^2w-(\alpha_1+\alpha_2)x-\beta_2z,\\
&z_{23}=z(zw+\beta_2)(-xzw+z^2w-\alpha_4x-\beta_2x+\beta_2z+\beta_4x)/x,\\
&w_{23}=-1/(z(zw+\beta_2)),
\end{split}\\
\begin{split}
&r_1(r_0r_3):\\
&x_{24}=1/((y+w)(tw-xw+ty-xy+\alpha_0)),\\
&y_{24}=(y+w)((t-x)(y+w)+\alpha_0)(tyw-xyw+ty^2-xy^2+\alpha_0y-(\alpha_1+\alpha_2)w)/w,\\
&z_{24}=tw^2+2(t-x)yw+ty^2-xy^2-zw^2+\alpha_0(y+w)+(\alpha_4-\beta_4)w, \quad w_{24}=1/w.
\end{split}
\end{align*}
\end{proposition}
Each of them contains a two-parameter family of meromorphic solutions of \eqref{1} as the initial conditions. 

It is still an open question whether the phase space of the system \eqref{1} can be covered by the original coordinate system $(x,y,z,w)$ in addition to the canonical coordinate systems $(x_i,y_i,z_i,w_i) \ (i=0,1, \ldots ,24)$.

\end{document}